%% file: limhydFHNCFFv1.tex
\documentclass[a4paper,11pt]{article}
\usepackage{amsmath}
\usepackage{amssymb}
\usepackage{setspace}
\usepackage{graphicx}
\usepackage{psfrag}
\usepackage{listings} 

\usepackage{ucs}
\usepackage[utf8x]{inputenc}
\usepackage{graphicx}
\usepackage{amsfonts}
\usepackage{dsfont}
\usepackage{amssymb}
\usepackage{amsmath}
\usepackage{amsthm}
\usepackage{enumerate}
\usepackage{stmaryrd}
\usepackage{fullpage}
\usepackage{ifthen}
\usepackage{subfigure}
\usepackage{epic}
\usepackage{authblk}
\usepackage{textcomp}
\usepackage[all]{xy}
\usepackage{mathrsfs}
\usepackage{color}
\usepackage{titlesec}

\usepackage[hypertexnames=false,colorlinks=true,linkcolor=blue,citecolor=blue]{hyperref}
\usepackage[numbers,comma,square,sort&compress]{natbib}

\newcommand{\vertiii}[1]{{\left\vert\kern-0.25ex\left\vert\kern-0.25ex\left\vert
          #1
        \right\vert\kern-0.25ex\right\vert\kern-0.25ex\right\vert}}

\newtheorem{lem}{Lemma}[section]
\newtheorem{thm}[lem]{Theorem}
\newtheorem{prop}[lem]{Proposition}
\newtheorem{cor}[lem]{Corollary}
\newtheorem{rmk}[lem]{Remark}
\newtheorem{defi}[lem]{Definition}

\newcommand{\etal}{\textit{et al.}\ }
\newcommand{\Supp}{\mathop{\mathrm{Supp}}}

\newcommand{\R}{\mathbb{R}}
\newcommand{\md}{\mathrm{d}}
\newcommand{\K}{\mathcal{K}}
\newcommand{\mD}{\mathcal{D}}
\newcommand{\mZ}{\mathcal{Z}}
\newcommand{\mS}{\mathcal{S}}
\newcommand{\mR}{\mathcal{R}}
\newcommand{\mF}{\mathcal{F}}
\newcommand{\mI}{\mathcal{I}}
\newcommand{\ve}{\varepsilon}
\newcommand{\dv}{\mathop{\mathrm{div}}}
\newcommand{\mx}{\mathbf{x}}

  {\begin{trivlist}\item[]{\bf Hypothesis #1 }\em}{\end{trivlist}}

\makeatletter\@addtoreset{figure}{section}\makeatother

\makeatletter \@addtoreset{equation}{section} \makeatother

\newcommand{\ds}{\displaystyle}

\newenvironment{Proof}[1][.]%
 {\begin{trivlist}\item[]\textbf{Proof#1 }}%
 {\hspace*{\fill}$\rule{0.3\baselineskip}{0.35\baselineskip}$\end{trivlist}}

\title{Rigorous derivation of the nonlocal reaction-diffusion FitzHugh-Nagumo system}
\author[1]{Joachim Crevat}
\author[1]{Gr\'egory Faye\footnote{Corresponding author: \texttt{gregory.faye@math.univ-toulouse.fr}}}
\author[1,2]{Francis Filbet}
\affil[1]{Institut de Math\'ematiques de Toulouse ; UMR5219, Universit\'e de Toulouse ; UPS IMT, F-31062 Toulouse Cedex 9 France}
\affil[2]{Institut Universitaire de France}

\begin{document}
\maketitle

\begin{abstract}
We introduce a spatially extended transport kinetic FitzHugh-Nagumo model with
forced local interactions and prove that its hydrodynamic limit
converges towards the classical nonlocal reaction-diffusion
FitzHugh-Nagumo system. Our approach is based on a relative entropy
method, where the macroscopic quantities of the kinetic model are
compared with the solution to the  nonlocal reaction-diffusion
system. This approach allows to make the rigorous link between
kinetic and reaction-diffusion models.
\end{abstract}


\section{Introduction}
\input{introduction.tex}

\section{Preliminaries and main result}\label{sec:main}
\input{mainresult.tex}

\section{A priori estimates}\label{sec:apriori}
\input{apriori.tex}

\section{Relative entropy estimate and proof of Theorem \ref{thm:2}}\label{sec:proof}
\input{entropy.tex}

\appendix

\section{Proof of Proposition \ref{prop:exfeps}}\label{appexfeps}

\input{AppThm1.tex}

\section{Proof of Proposition \ref{prop:exhydro}}\label{appexhydro}
\input{AppThm2.tex}

\bibliography{plain}

\end{document}

%% file: introduction.tex
Generally, neuron models focus on the regulation of the electrical potential of the membrane of a nerve cell depending on the input it receives. This regulation is the result of ionic exchanges between the neuron and its environment through its cellular membranes. A very precise modeling of these ion exchanges led to the well-known Hodgkin-Huxley model \cite{HOD}. In this paper, we shall rather focus on a simplified version, called the FitzHugh-Nagumo (FHN) model \cite{FIT,NAG}, which keeps its most valuable aspects, but remains relatively simple mathematically. More precisely, the FHN model accounts for the variations of the membrane potential $v$ of a neuron coupled to an auxiliary variable $w$ called the adaptation variable. It is usually written as follows
\begin{equation}
\left\{
\begin{array}{l}
 \ds\dfrac{\md v}{\md t} \,=\,N(v) \,-\, w \,+\, I_{\text{ext}},\\ \, \\
\ds\dfrac{\md w}{\md t} \,=\,  \tau\,(v\,+\,a \,-\, b\,w),
\end{array}\right.
\label{eqFHN}
\end{equation}
where $I_{\text{ext}}$ stands for the input current the neuron receives, $N(v)$ is a bistable nonlinearity which models the cell excitability, and $\tau\geq0$, $a\in\R$ and $b\geq0$ are some given constants. Without loss of generality \cite{BAL,MIS,NAG,FAY}, we assume that $N(v)$ is given by the following cubic nonlinearity\footnote{All our results remain true for $N(v)=v(\alpha-\beta|v|^2)$ with any $\alpha$, $\beta>0$.}
\begin{equation}
N(v):=v-v^3, \quad v\in\R.
\label{cubic}
\end{equation}
Another mathematical reason for looking at system \eqref{eqFHN} is that it is a prototypical model of excitable kinetics and interest in such systems stems from the fact that although the kinetics are very simple spatial coupling can produce complex dynamics, where well-known examples are the propagation of excitatory pulses, spiral waves in two-dimensions, and spatio-temporal chaos. Here, we introduce spatial coupling through the input current $I_{\text{ext}}$ which models how the membrane potentials of the other neurons influence the membrane potential of one given neuron. More specifically, following \cite{BAL,MIS}, we consider a network of $n$ neurons interacting through their synapses, which we suppose to be of electrical type, such that the evolution of the pair voltage-adaptation $(v_i,w_i)_{1\leq i \leq n}$ satisfies the system of equations
\begin{equation}
\left\{
\begin{array}{l}
 \ds\dfrac{\md v_i}{\md t} \,=\,N(v_i) \,-\, w_i \,-\, \dfrac{1}{n}\overset{n}{\underset{j=1}{\sum}} \Phi_\ve(\mathbf{x}_i-\mathbf{x}_j) (v_i-v_j),\\ \,\\
 \ds\dfrac{\md w_i}{\md t} \,=\,  \tau(v_i \,+\, a \,-\, b w_i),
\end{array}\right.
\label{eqsyst}
\end{equation}
where we choose to associate to each neuron $i \in \{1,...,n\}$ its position in the network $\mathbf{x}_i\in\mathbb{R}^d$ with $d\in\{1,2,3\}$ as in \cite{LUC}. The connectivity kernel $\Phi_\ve:\mathbb{R}^d\rightarrow\mathbb{R}$ models the effect of the distance between two neurons on their electro-chemical interactions, and here $\ve>0$ is a scaling parameter as explained later. System \eqref{eqsyst} is a natural spatial extension of the FHN models introduced and studied in \cite{BAL,MIS}.

The main purpose of this article is to establish a macroscopic description of the neural network \eqref{eqsyst} in the limit $n\rightarrow +\infty$, choosing an appropriate scaling for the connectivity function $\Phi_\ve $. Our strategy is to introduce an intermediate modeling scale, namely a mesoscopic description of system \eqref{eqsyst}. We thus associate to system \eqref{eqsyst} the following spatially-extended kinetic equation
\begin{equation}\label{eq:f}
 \partial_t f^\ve  + \partial_v \left (f^\ve  N(v) - f^\ve w - f^\ve \K_{\Phi_\ve}[f^\ve]  \right) + \partial_w \left( A(v,w)f^\ve \right) = 0, \quad \text{ on } (0,\infty)\times \R^{d+2},
\end{equation}
where the operator $\K_{\Phi_\ve}[f^\ve]$ and the convective term $A$ are given by  
$$
\left\{
\begin{array}{l}
\ds\K_{\Phi_\ve}[f^\ve](t,\mathbf{x},v) \,:=\, \int_{\mathbb{R}^{d+2}} \Phi_\ve(\mathbf{x}-\mathbf{x}')(v-v')\, f^\ve(t,\mathbf{x}',v',w')\,\md w' \md v' \md \mathbf{x}',
\\ \, \\
 \ds A(v,w) \,:=\, \tau(v \,+\, a \,-\, b\,w).
\end{array}\right.
$$
Such an equation models the evolution of the density function $f^\ve(t,\mathbf{x},v,w)$ of finding neurons with a potential membrane $v\in\mathbb{R}$ and an adaptation variable $w\in\mathbb{R}$ at time $t\geq0$ and position $\mathbf{x} \in \R^d$ within the cortex. The term $\K_\ve[f^\ve]$ describes nonlocal interactions through the whole field, whereas the other terms account for the local reactions due to the excitability of nerve cells. 

Using classical arguments \cite{CAN,GOL} and under some mild assumptions on $\Phi_\ve$, it can actually be proved that the kinetic equation (\ref{eq:f}) is the mean-field limit of (\ref{eqsyst}) as $n$ goes to infinity \cite{CRE:18}. In fact, the mean-field limit of the  FHN system \eqref{eqsyst} with stochastic fluctuations towards a kinetic PDE similar to (\ref{eq:f}) has already been studied for conductance-based connectivity without spatial dependance \cite{BAL,BFT}, for homogeneous interactions \cite{MIS} and  for spatially organized neural network with compactly supported communication weight \cite{LUC}. Various other types of kinetic models have been derived during the past decades depending on the hypotheses assumed for the dynamics of the emission of an action potential. They include for example integrate-and-fire neural networks \cite{brunel,CCP,CPSS} and time-elapsed neuronal models \cite{PPS,CCDR,CHE1,CHE}. 

Coming back to our problem of establishing a macroscopic description of the neural network \eqref{eqsyst}, we introduce the following macroscopic quantities for $f^\ve$ a solution of \eqref{eq:f}:
\begin{equation}
\left\{
\begin{array}{l}
\ds\rho^\ve(t,\mathbf{x}) \,:=\,\displaystyle\int_{\R^2}f^\ve(t,\mathbf{x},v,w)\,\md v \md w,   \\ \, \\
\ds\rho^\ve(t,\mathbf{x})V^\ve(t,\mathbf{x}) \,=\, j^\ve(t,\mathbf{x})\,:=\;\displaystyle\int_{\R^2}vf^\ve(t,\mathbf{x},v,w)\,\md v \md w,   \\ \, \\
\ds\rho^\ve(t,\mathbf{x})W^\ve(t,\mathbf{x}) \,=\, q^\ve(t,\mathbf{x})\,:=\;\displaystyle\int_{\R^2}wf^\ve(t,\mathbf{x},v,w)\, \md v\md w,   
\end{array}
\label{eqMacro}
\right.
\end{equation}
where $\rho^\ve$ stands for the averaged neuron density, $V^\ve$ is the average membrane potential, and $W^\ve$ is the average adaptation variable.  Formally, it is not difficult to check that the equation satisfied by $j^\ve$ is not closed because of the presence of the reaction term $\partial_v \left( f^\ve(N(v)\right)$ in \eqref{eq:f}, which introduces higher moments of $f^\ve$ in $v$. To circumvent this difficulty and obtain a closed macroscopic system, we consider a specific form for the connectivity function $\Phi_\ve$. We will assume that it can be decomposed by the superposition of strong local interactions modeled by $\frac{1}{\ve}\delta_0$, where $\delta_0$ is the Dirac distribution centered at $0$ and weak lateral interactions described by a non-negative connectivity kernel $\Psi$. As a consequence, for all $\ve>0$, we work with 
\begin{equation}
\Phi_\ve=\Psi+\frac{1}{\ve}\delta_0.
\label{eqPhieps}
\end{equation}
Note that such an assumption of strong local interactions and weak lateral interactions is often used when modeling visual cortex \cite{BRE03}. We can then rewrite the kinetic equation \eqref{eq:f} as
\begin{equation}
\label{eq:feps}
 \partial_t f^{\varepsilon}  + \partial_v \left (f^{\varepsilon}  N(v) - f^{\varepsilon} w - f^{\varepsilon} \K_{\Psi}[f^{\varepsilon}]  \right) - \displaystyle\frac{1}{\varepsilon}\partial_v \left(f^{\varepsilon}( \rho^\ve v-j^{\varepsilon} ) \right) + \partial_w \left( A(v,w)f^{\varepsilon} \right) = 0
\end{equation}
and the macroscopic quantities $(\rho^{\varepsilon},\rho^\ve V^{\varepsilon},\rho^\ve W^{\varepsilon})$ formally verify the system of equations:
\begin{equation}
\left\{
\begin{array}{l}
\ds\partial_t \rho^{\varepsilon} \,=\, 0, \\ \, \\
\ds\partial_t\left(\rho^{\varepsilon}V^{\varepsilon}\right)\,-\,\rho^{\varepsilon}\mathcal{L}_{\rho^{\varepsilon}}(V^{\varepsilon})  \,=\,\rho^{\varepsilon}\left[ N(V^\varepsilon)\,-\, W^{\varepsilon}\right]  \,+\,\mathcal{E}(f^\varepsilon),  \\ \, \\
\ds\partial_t\left(\rho^{\varepsilon}W^{\varepsilon}\right) \,=\, \rho^{\varepsilon}\,A(V^{\varepsilon},W^{\varepsilon}),
\end{array}\right.
\label{eq:feps2}
\end{equation}
where the spatial nonlocal operator $\mathcal{L}_\rho$ is defined through
\begin{equation}
\mathcal{L}_\rho(V)\, :=\, - (\Psi\ast\rho)V \,+\, \Psi\ast[\rho V],
\label{opL}
\end{equation}
whereas $\ast$ stands for the convolution product in $\mathbb{R}^d$ only with respect to space. Furthermore, we have set
\begin{equation}
\mathcal{E}(f^\varepsilon)(t,\mathbf{x}):= \int_{\mathbb{R}^2} f^{\varepsilon}(t,\mathbf{x},v,w)\left[ N(v)\,-\,N(V^\varepsilon(t,\mathbf{x}))\right]\, \md v\md w.
\label{errorR}
\end{equation}
Let us formally determine the hydrodynamic limit of the equation (\ref{eq:feps}) as $\varepsilon\rightarrow0$.  If we multiply (\ref{eq:feps}) by $\varepsilon$ and if we take $\varepsilon\rightarrow0$, the local interaction term $\rho^\ve (v-V^{\varepsilon})f^{\varepsilon}$ vanishes. Consequently,  $f^{\varepsilon}$ should converge in some weak sense to the mono-kinetic distribution in $v$
\begin{equation}
f^{\varepsilon}(t,\mathbf{x},v,w) \underset{\varepsilon \rightarrow 0}{\rightharpoonup} F(t,\mathbf{x},w) \otimes \delta_{0}(v-V(t,\mathbf{x})),
\label{weaklimit}
\end{equation}
 where the limit $(F,V)$ has to be determined. First we set
$$
\rho(t,\mathbf{x}) \,:=\, \int_{\mathbb{R}} F(t,\mathbf{x},w) \md w, \quad \rho(t,\mx)\, W(t,\mathbf{x}) \,:=\, \int_{\mathbb{R}} F(t,\mathbf{x},w) \,w \,\md w, 
$$
hence the macroscopic equation satisfied by  $V$ is formally obtained by observing that $\mathcal{E}(f^\varepsilon)\rightarrow 0$ as $\varepsilon\rightarrow0$ in the second equation of  (\ref{eq:feps2}), hence we get
$$
\partial_t\left(\rho \,V\right)\,-\,\rho\mathcal{L}_{\rho}(V)  \,=\,\rho\,\left[ N(V)\,-\, W\right]. 
$$
Then substituting (\ref{weaklimit}) in the kinetic equation  (\ref{eq:feps}) and integrating with respect to $v\in\mathbb{R}$, it yields that $F$ is solution to
\begin{equation*}
\partial_t F + \partial_w \left( A(V,w) \, F \right) \,=\, 0.
\end{equation*}
Therefore, the limit function $\mZ:=(\rho,\rho V,\rho W)$ is expected to verify the following nonlocal reaction diffusion system:
\begin{equation}
\label{eq:lim}
\partial_t \mZ = \mF(\mZ), 
\end{equation}
with $\mF$ given by
\begin{equation}
\label{def:F}
\mF(\mZ) \,:=\,   \rho\,\begin{pmatrix} 0 \\ \mathcal{L}_\rho(V) \,+\,
  N(V) \,-\, W  \\ A(V,W) \end{pmatrix}.
\end{equation}

Let us first make some  comments on the structure of the macroscopic model \eqref{eq:lim}-(\ref{def:F}). For all $\mx\in\R^d$ such that $\rho(t,\mathbf{x})=\rho_0(\mathbf{x})>0$, the system  \eqref{eq:lim} reduces to the usual nonlocal reaction-diffusion system of FitzHugh-Nagumo type
\begin{equation}
\left\{
\begin{array}{l}
 \partial_t V \,-\,\mathcal{L}_{\rho_0}(V)  \,=\,  N(V) \,-\,  W,\,\\ \, \\
 \partial_t W \,=\, \tau(V+a-bW), 
\end{array}\right.
\label{eq:hydro}
\end{equation}
where $\mathcal{L}_{\rho_0}(V)$ can be interpreted as a nonlocal diffusion operator in $\mathbf{x}$. In the limiting case $\rho_0 \equiv 1$\footnote{In our setting, $\rho_0$ is probability density function such that the case $\rho_0 \equiv 1$ is excluded from our hypotheses. Nevertheless, system \eqref{eq:hydro} is still well defined for $\rho_0 \equiv 1$, see Proposition \ref{prop:exhydro}.}, such a system has already been well studied especially regarding the formation of propagating waves (traveling fronts and pulses) in both cases $\tau=0$ and  $0<\tau\ll 1$ \cite{bates-etal:97,FAY}. We also mention the classical works regarding the local FitzHugh-Nagumo system, that is when $\mathcal{L}_{\rho_0}(V)$ is replaced by the standard diffusion operator \cite{CAR,HAS,JON}, and the more recent advances for the discrete case \cite{HUP10,HUP13}. The present work is then a rigorous justification of the nonlocal reaction-diffusion system of FitzHugh-Nagumo type \eqref{eq:lim} that is obtained as the hydrodynamic limit of the kinetic equation \eqref{eq:feps} as $\varepsilon\rightarrow0$.

The main challenge towards a rigorous proof of this hydrodynamic limit stems from the Dirac singularity of the mono-kinetic distribution which prevents us from using a classical entropy of the form $f\log(f)$ since it would not be well-defined. Following ideas from \cite{FIG}, who proved a similar hydrodynamic limit of the kinetic Cucker-Smale model for collective motion with forced local alignment towards the pressureless Euler equation with a nonlocal force, we shall overcome this problem by the mean of a relative entropy argument. Let us also mention the work of \cite{KAN} in which the hydrodynamic limit of a collisionless and non-diffusive kinetic equation under strong local alignment regime is rigorously established \textit{via} a relative entropy argument. The specific difficulty here is that, instead of having a transport term as in \cite{FIG}, the presence of the reaction term $\partial_v(f^\varepsilon N(v))$ in \eqref{eq:feps} introduces higher order moments of $f^\varepsilon$ in $v$ which we will need to control. In fact, it will be enough to have a priori estimates of second and fourth order moments  of $f^\ve$ in $(v,w)$ to circumvent this difficulty.

We conclude this introduction by comparing the kinetic model  \eqref{eq:f} to the pioneer works of Amari, Wilson and Cowan in the 1970s \cite{AMA,WC1,WC2} who heuristically derived spatially extended models describing the macroscopic activity of large assemblies of neurons. Such models are referred to as neural field equations in the literature and account for the time evolution of the spatial averaged membrane potential of a population of neurons.  Such a class of models has received much attention and has been very successful at reproducing a number of biological phenomena, including in particular visual hallucinations, binocular rivalry or working memory. We refer to the recent surveys \cite{bressloff:12, bressloff:14,coombes:05} for more developments on neural field models and applications in neuroscience. By derivation, neural field equations are macroscopic models and it is still an open problem to rigorously justify such models from conductance-based neural networks. In a slightly different direction, we mention the recent work of Chevallier \etal \cite{CHE} on mean-field limits for nonlinear spatially extended Hawkes processes with exponential memory kernels who are able to recover at the limit a process whose law is an inhomogeneous Poisson process having an intensity which solves a scalar neural field equation. In this article, we have addressed this open question in the case of the FHN model.

The rest of this article is organized as follows. In Section \ref{sec:main}, we state our main result about the hydrodynamic limit of (\ref{eq:feps}) towards (\ref{eq:lim}). Then, in Section \ref{sec:apriori}, we derive some a priori estimates that will be crucial for our relative entropy argument. The proof of our main result is contained in Section \ref{sec:proof}.

%% file: mainresult.tex
In this section, we present our main result on the hydrodynamic limit
from a weak solution $(f^\ve)_{\ve >0}$ of the kinetic equation \eqref{eq:feps} to a classical
solution $(\rho_0,\rho_0 V,\rho_0 W)$ of the asymptotic system (\ref{eq:lim}). For that, we first
need to present the existence result for the weak solution of \eqref{eq:feps}
and the classical solution of (\ref{eq:lim}).

First we set the hypotheses we make for the study of the kinetic FHN
equation (\ref{eq:feps}) and the limit system (\ref{eq:lim}). We consider a connectivity kernel $\Psi$ in \eqref{eqPhieps} which  is
non-negative, symmetric and satisfies
\begin{equation}
\label{H1}
\Psi \in L^1(\R^d).
\end{equation}
This last assumption models the fact that if two neurons are far away from each other, they have weak mutual interactions. The other condition is a natural biological assumption and expresses the symmetric and excitatory nature of the considered underlying neural network.

\subsection{Existence of weak solution to (\ref{eq:feps})}
We here say that $f^\ve$ is a weak solution of (\ref{eq:feps}) if for
any $T >0$, $f^\ve(0,.) = f^\ve_0\geq 0$  in  $\mathbb{R}^{d+2}$,
$$
f^\ve \in \mathscr{C}^0([0,T], L^1(\mathbb{R}^{d+2})) \cap L^\infty((0,T)\times\mathbb{R}^{d+2}),
$$
and (\ref{eq:feps}) holds in the sense of distribution, that is, for any $\varphi\in C_c^\infty([0,T)\times\mathbb{R}^{d+2})$, the weak formulation holds
\begin{eqnarray}
\nonumber
&&\int_0^T \int  f^{\varepsilon} \,\left[\partial_t\varphi + \left (N(v)
    - w - \K_{\Psi}[f^{\varepsilon}]  -\frac{1}{\varepsilon}( \rho^\ve
    v-j^{\varepsilon} ) \right)\partial_v\varphi  + A(v,w)\partial_w
  \varphi \right]\md \mathbf{z}\,\md t 
\\
&&+\, \int  f^{\varepsilon}_0 \varphi(0)\,\md \mathbf{z} \,=\, 0
\label{weak:sol}
\end{eqnarray}
where $\mathbf{z}=(\mathbf{x},v,w)\in\mathbb{R}^{d+2}$.

\begin{rmk}
Using the mass conservation property of equation \eqref{eq:feps}, we can easily check that the time varying macroscopic quantities $(\rho^\ve,V^\ve,W^\ve)$ defined in \eqref{eqMacro} simplify to
\begin{equation}
\rho^{\varepsilon}(t,\mathbf{x})\,:=\,\ds\int_{\mathbb{R}^2}f^{\varepsilon}(t,\mathbf{x},v,w)\,\md v\, \md w
=\int_{\mathbb{R}^2}f^{\varepsilon}_0(\mathbf{x},v,w)\,\md v\, \md
w=\rho_0^\ve(\mathbf{x}),
\label{rho:constant}
\end{equation}
hence we have
$$
\left\{\begin{array}{l}
\rho_0^{\varepsilon}(\mathbf{x})V^{\varepsilon}(t,\mathbf{x}) \,:=\, \ds\int_{\mathbb{R}^2} v \, f^{\varepsilon}(t,\mathbf{x},v,w)\, \md v\,\md w,\\ \, \\
\rho_0^{\varepsilon}(\mathbf{x})W^{\varepsilon}(t,\mathbf{x}) \,:=\; \ds\int_{\mathbb{R}^2} w \, f^{\varepsilon}(t,\mathbf{x},v,w)\, \md v\,\md w,
\end{array}\right.
$$
for all $\mathbf{x}\in\R^d$ and all $t >0 $ where $f^\ve$ is well-defined. 
\end{rmk}

First, let us prove the
well-posedness of the kinetic equation \eqref{eq:feps}. 

\begin{prop}\label{prop:exfeps}
For any $\ve>0$ we choose $\Psi$ to be non-negative, symmetric and satisfies
(\ref{H1}), we also  assume that $f^\ve_0$ satisfies
\begin{equation}
\label{H3-1}
f^\ve_0 \geq 0, \quad f^\ve_0\in L^1(\mathbb{R}^{d+2}), \quad f^\ve_0,\, \nabla_{\mathbf{u}} f_0^\ve \in L^\infty(\mathbb{R}^{d+2}),
\end{equation}
where $\mathbf{u}=(v,w)$ and for all $\mathbf{x}\in\R^d$, 
\begin{equation}
\label{H3-2}
\Supp(f_0^\ve(\mathbf{x},\cdot,\cdot))\subseteq B(0,R_0^\ve)\subset \R^2.
\end{equation}
Then for any $T>0$, there exists a unique  $f^\ve$ weak solution to
\eqref{eq:feps} in the sense of (\ref{weak:sol}), which is compactly
supported in $\mathbf{u}=(v,w)\in\mathbb{R}^2$.
\end{prop}

The proof relies on a classical fixed point argument, but for the sake of completeness, it is postponed to the Appendix \ref{appexfeps}.

\subsection{Existence of classical solution to the nonlocal FitzHugh-Nagumo system}
Let us now state the result of existence and uniqueness for the nonlocal reaction-diffusion FitzHugh-Nagumo system defined as

\begin{equation}
\left\{
\begin{array}{ll}
 \partial_t V \,-\,\mathcal{L}_{\rho_0}(V)  \,=\,  N(V) \,-\,  W, \,\\ & t>0 \text{ and } \mx \in\R^d, \\
 \partial_t W \,=\, \tau(V+a-bW), & \,\\ \, \\
 V(0,\mx)=V_0(\mx), \quad W(0,\mx)=W_0(\mx), & \mx \in\R^d.
\end{array}\right.
\label{eq:nlFHN}
\end{equation}

Before describing precisely the  existence  and uniqueness of
solution  $(V,W)$ to the hydrodynamical system (\ref{eq:nlFHN}), let
us emphasize that this system is more convenient to analyse than \eqref{eq:lim}-\eqref{def:F}
verified by $\mZ=(\rho_0, \,\rho_0\,V,\,\rho_0\,W)$. Indeed, as we will
see, both solutions coincide in the region of interest where $\rho_0>0$,
but the study of (\ref{eq:nlFHN}) allows to construct a
solution such that  for all $t\in[0,T]$
$$
V(t), \, W(t) \in L^\infty(\R^d).
$$ 
This property is crucial to apply the relative entropy method
in the asymptotic analysis of \eqref{eq:feps} when $\ve\rightarrow 0$.

\begin{prop}\label{prop:exhydro}
We choose $\Psi$ to be non-negative, symmetric and satisfies
(\ref{H1}), we also suppose that $\rho_0$ and the initial data $(V_0,W_0)$
satisfies, 
\begin{equation}
\label{H4-1}
\rho_0\geq 0,\quad\rho_0 \in L^1\cap L^\infty (\mathbb{R}^{d}),\quad
 V_0,W_0\in L^{\infty}(\mathbb{R}^d).
\end{equation}
 Then for any $T>0$, there exists a unique classical solution $(V,W)\in
 \mathscr{C}^1([0,T], L^{\infty}(\mathbb{R}^d))$ to  the equation
 \eqref{eq:nlFHN}. Furthermore, $\mathcal{Z}=(\rho_0,\rho_0 V,\rho_0 W)$
 is a solution to \eqref{eq:lim}-(\ref{def:F}).
 \end{prop}

The proof of this proposition also relies on a classical fixed point
argument. For the sake of completeness, it is postponed to the
Appendix \ref{appexhydro}. Then as a direct consequence, for a given
initial data $(F_0,V_0)$, it leads to the existence and uniqueness of solution $(F,V)$ to
the following system of equations
\begin{equation}
\label{eq:F}
\left\{
\begin{array}{l}
\ds \partial_t F + \partial_w \left( A(V,w) \, F \right) \,=\, 0,
\\ \, \\
\rho_0\ds \partial_tV\,-\,\rho_0\mathcal{L}_{\rho_0}(V)  \,=\,\rho_0 \left[ N(V)\,-\, W\right],
\\ \, \\
\rho_0(\mathbf{x}) = \ds \int_{\mathbb{R}} F(t,\mathbf{x},w)\,\md w,
\\ \, \\
V(t,\mx)\, = \,0, \quad t>0 \text{ and } \mx\in \R^d \, \backslash \, \text{Supp}(\rho_0),
\\ \, \\
W(t,\mathbf{x}) = \left\{
\begin{array}{ll}
\ds\frac{1}{\rho_0(\mathbf{x})}\int_{\mathbb{R}} F(t,\mathbf{x},w) \,w\,\md w, &
                                                               \textrm{if }\, \rho_0(\mathbf{x})>0,
\\ \, \\
0, &
                                                               \textrm{else.}
\end{array}\right.
\end{array}\right.
\end{equation} 

More precisely we have the following result.
 
\begin{cor}
\label{cor:exF}
Consider $V_0\in L^\infty(\R^d)$ and $F_0\in\mathcal{M}(\mathbb{R}^{d+1})$ such that 
$$
\int_{\mathbb{R}^{d+1}} w^2 F_0(\md {\bf x},\md w) < \infty
$$
and for almost every $\mathbf{x}\in\mathbb{R}^d$,
$$
\rho_0(\mathbf{x})\,=\, \int_{\mathbb{R}} F_0(\mathbf{x},w) \,\md w, \quad W_0(\mathbf{x}) := \left\{
\begin{array}{ll}
\ds\frac{1}{\rho_0(\mathbf{x})}\int_{\mathbb{R}} F_0(\mathbf{x},w) \,w\,\md w, &
                                                               \textrm{if }\, \rho_0(\mathbf{x})>0,
\\ \, \\
0,  &
                                                               \textrm{else,}
\end{array}\right.
$$
where $\rho_0 \in L^1\cap L^\infty (\mathbb{R}^{d})$ and $W_0\in L^\infty(\R^d)$. We further assume that $V_0=0$ on $ \R^d \, \backslash \, \mathrm{Supp}_\text{ess}(\rho_0)$.

Then for any $T>0$, there exists a unique couple $(F,V)$ solution to \eqref{eq:F}, where
$$
(F,V) \in
L^\infty((0,T), \mathcal{M}(\mathbb{R}^{d+1}))\times
\mathscr{C}^1([0,T], L^{\infty}(\mathbb{R}^d)),
$$ 
and $F$ is a measure solution to the first equation 
\eqref{eq:F}, that is, for any $\varphi\in\mathscr{C}^1_c(\mathbb{R}^{d+1})$
\begin{equation}
\label{meas:0}
\frac{\md}{\md t}\int_{\mathbb{R}^{d+1}} \varphi(\mathbf{x},w) \, F(t,\md \mathbf{x},\md
w) \,-\, \int_{\mathbb{R}^{d+1}} 
A(V(t,\mathbf{x}),w) \,\partial_w\varphi(\mathbf{x},w) \, F(t,\md \mathbf{x},\md w)  \,=\,0  
\end{equation}
such that there exists a constant $C_T>0$, 
$$
\int_{\mathbb{R}^{d+1}} w^2 F(t,\md {\bf x},\md w) < C_T, \quad t\in [0,T].
$$ 
\end{cor}
\begin{Proof}
We first apply Proposition \ref{prop:exhydro} with $(\rho_0,V_0,W_0)$,
where $(\rho_0,W_0)$ is computed from the moments with respect to
$(1,w)$  of the initial distribution $F_0$, that is, for almost every $\mx\in\R^d$,
$$
\rho_0(\mathbf{x})\,=\, \int_{\mathbb{R}} F_0(\mathbf{x},w) \,\md w, \quad W_0(\mathbf{x}) := \left\{
\begin{array}{ll}
\ds\frac{1}{\rho_0(\mathbf{x})}\int_{\mathbb{R}} F_0(\mathbf{x},w) \,w\,\md w, &
                                                               \textrm{if }\, \rho_0(\mathbf{x})>0,
\\ \, \\
0,  &
                                                               \textrm{else.}
\end{array}\right.
$$
We then denote by $(\widetilde{V},\widetilde{W})$ the corresponding unique classical solution to (\ref{eq:nlFHN}) starting from such an initial condition. 
Finally, we define
\begin{equation}
\label{eq:Vnew}
V(t,\mx) := \left\{
\begin{array}{ll}
\widetilde{V}(t,\mx), &
                                                               \textrm{if }\, \rho_0(\mathbf{x})>0,
\\ \, \\
0, &
                                                               \textrm{else.}
\end{array}\right.
\end{equation}
As $\rho_0$ is independent of time, we note that $V\in\mathscr{C}^1([0,T], L^{\infty}(\mathbb{R}^d))$.

We now prove that for any couple solution $(\overline{F}, \overline{V})$ of \eqref{eq:F} then $\overline{V}$ is precisely given by \eqref{eq:Vnew}. Thus, let us suppose that $(\overline{F}, \overline{V})$ is a well-defined solution of \eqref{eq:F} on $[0,T]$  with finite second moment. We necessarily get that $\rho_0\overline{W}$ has to satisfy 
$$
\partial_t (\rho_0 \overline{W}) -  \rho_0 \,A(\overline{V},\overline{W})  \,=\, 0,
$$
together with
$$
\rho_0\ds \partial_t \overline{V}\,-\,\rho_0\mathcal{L}_{\rho_0}(\overline{V})  \,=\,\rho_0\left[ N(\overline{V})\,-\, \overline{W}\right].
$$
Thus, for any $\mx\in\R^d$ such that $\rho_0(\mx)>0$, the couple $(\overline{V},\overline{W})$ coincides with the unique solution to \eqref{eq:nlFHN} with initial condition $(\rho_0,V_0,W_0)$. This is ensured from the fact that the convolution in the nonlocal part $\Psi*(\rho_0 \overline{V})$ of the linear operator $\mathcal{L}_{\rho_0}(\overline{V})$ is only evaluated on the regions where $\rho_0>0$. Then, this uniquely defines $\overline{V}(t,\mx) = \widetilde{V}(t,\mx)$ on $\rho_0>0$. By definition of a solution to \eqref{eq:F}, whenever $\rho_0=0$, we have that $\overline{V}(,\mx)=0$.
As a conclusion, we have just shown that if $(\overline{F},\overline{V})$ is a well-defined solution of \eqref{eq:F} on $[0,T]$  with finite second moment then necessarily $\overline{V} = V$ where $V$ is uniquely defined in \eqref{eq:Vnew}.

For $V$ defined as above in \eqref{eq:Vnew}, we consider the transport equation
\begin{equation}
\label{eq:transportF}
\left\{
\begin{array}{l}
\ds \partial_t F + \partial_w \left( A(V,w) \, F \right) \,=\, 0,
\\ \, \\
F(0)  \,=\,F_0.
\end{array}\right.
\end{equation} 
Then, for almost every $\mathbf{x}\in\mathbb{R}^d$ and all $(t,w)\in[0;T]\times\R$, we introduce 
the system of characteristic curves  associated to 
\eqref{eq:transportF}, that is,
\begin{equation}
\label{eq:characF}
\left\{
\begin{array}{l}
\ds\frac{\md}{\md s} \mathcal{W}(s) \,=\, A\left(V(s,\mathbf{x})\, ,\,\mathcal{W}(s)  \right),   \\ \, \\
\mathcal{W}(t)\,=\,w,
\end{array}
\right.
\end{equation}
where $V$ is defined in \eqref{eq:Vnew}. From the regularity with respect to
$(t,w)\in [0,T]\times\mathbb{R}$ of the functions $A(.,.)$ and $V(.,.)$ and since $A(.,.)$ grows at most
linearly with respect to $w$, we get global existence and uniqueness of a solution to
\eqref{eq:characF}. This solution is denoted by
$\mathcal{W}(s,t,\mathbf{x},w)$, then we verify using the theory of characteristics
that the unique solution to the transport equation  \eqref{eq:transportF} is given by
$$
F(t,\mathbf{x},w)= F_0(\mathbf{x},\mathcal{W}(0,t,\mathbf{x},w)) \,e^{\tau\, b\, t}.
$$
From the above expression, it easy to compute the second order moment of $F$ with respect to $w$
and get that there exists $C_T>0$ such that
$$
\int_{\mathbb{R}^{d+1}} w^2 F(t,\md {\bf x},\md w) < C_T, \quad t\in [0,T]. 
$$
As a final consequence of the above computations, we have that 
$$
\ds\frac{1}{\rho_0(\mathbf{x})}\int_{\mathbb{R}} F(t,\mathbf{x},w) \,w\,\md w = \widetilde{W}(t,\mx), \text{ for all } \mx\in\R^d \text{ such that } \rho_0(\mx)>0,
$$
by uniqueness of the solutions of \eqref{eq:nlFHN}. This shows that the couple $(F,V)$ is the unique solution to \eqref{eq:F} where $V$ is given in \eqref{eq:Vnew} and $F$ is the unique measure solution of \eqref{eq:transportF}.
\end{Proof}

\subsection{Main result}

Now, we are ready to state our main result about the hydrodynamic
limit. To this aim we  consider a non-negative initial data
$(f_0^{\varepsilon})_{\ve>0}$ and suppose that there exists a constant
$C>0$, such that for all $\ve>0$,

\begin{equation}
\label{H2-f0:L1}
  \|\rho_0^\ve\|_{L^\infty} \leq C
\end{equation}
and 
\begin{equation}
\label{H2-f0:M4}
 \int \left(1+|\mathbf{x}|^4 + |v|^4 + |w|^4\right) f_0^{\varepsilon}(\mathbf{x},v,w)\md w \md v \md \mathbf{x} \leq C.
\end{equation}
\begin{thm}
\label{thm:2}
Let $T>0$, $\Psi$ be a non-negative, symmetric kernel verifying (\ref{H1}) and
 $(f^\ve_0)_{\ve}$, a sequence of initial data such that for all $\ve>0$, \eqref{H3-1}, \eqref{H3-2}, \eqref{H2-f0:L1} and \eqref{H2-f0:M4} are satisfied.  Assume that  $(\rho_0,V_0, W_0)$ verifies \eqref{H4-1}, and
\begin{equation}
\label{H5}
\| \rho_0^\ve-\rho_0 \|_{L^2}^2  \,+\,
\int\rho_0^\ve(\mathbf{x})\left[ |V_0^{\varepsilon}(\mathbf{x}) - V_0(\mathbf{x}) |^2\,+\,
  |W_0^{\varepsilon}(\mathbf{x}) - W_0(\mathbf{x}) |^2 \right]\md \mathbf{x} \leq C \,\ve^{1/(d+6)}.
\end{equation}
Then, the macroscopic quantities $(\rho^\ve_0,V^\ve,W^\ve)$, computed from the solution
$f^\ve$ to \eqref{eq:feps}, verify that for all $t\in [0,T]$,
\begin{equation*}
\int_{\R^d} \left[ |V^{\varepsilon}(t,\mathbf{x}) -
  V(t,\mathbf{x})|^2 + |W^{\varepsilon}(t,\mathbf{x}) - W(t,\mathbf{x})|^2\right]\, \rho_0^{\varepsilon}(\mathbf{x})\,\md \mathbf{x} \,\leq C_T\, \ve^{1/(d+6)},
\end{equation*}
where $(V,W)$ is the unique solution to (\ref{eq:hydro}). 

Let further assume that $(\rho_0,V_0, W_0)$  are such that,
 $$
\rho_0(\mathbf{x})\,=\, \int_{\mathbb{R}} F_0(\mathbf{x},w) \,\md w, \quad W_0(\mathbf{x}) = \left\{
\begin{array}{ll}
\ds\frac{1}{\rho_0(\mathbf{x})}\int_{\mathbb{R}} F_0(\mathbf{x},w) \,w\,\md w, &
                                                               \textrm{if }\, \rho_0(\mathbf{x})>0,
\\ \, \\
0,  &
                                                               \textrm{else,}
\end{array}\right.
$$ 
for $F_0\in\mathcal{M}(\R^{d+1})$, and $V_0=0$ on $\R^d\backslash \Supp_\text{ess}(\rho_0)$. Moreover, consider the function $F_0^\ve$ such that for all $\mathbf{x}\in\R^d$ and all $w\in\R$,
$$F_0^\ve(\mathbf{x},w)\,=\,\ds\int f_0^\ve(\mathbf{x},v,w)\,\md v.$$
If  $F_0^\ve\rightharpoonup F_0$ weakly-$\star$ in $\mathcal{M}\left(\mathbb{R}^{d+1}\right)$
then we have for all $\varphi\in \mathscr{C}^0_b(\R^{d+2})$:
$$ 
\ds \int \varphi(\mathbf{x},v,w)\, f^{\varepsilon}(t,\mathbf{x},v,w)\,\md v\, \md w\, \md \mathbf{x} \rightarrow \int \varphi(\mathbf{x},V(t,\mathbf{x}),w)\, F(t,\md \mathbf{x},\md w)\, ,  
$$
strongly in
$L^1_{\text{loc}}(0,T)$ as $\varepsilon \rightarrow 0$, where $(F,V)$ is the unique solution to (\ref{eq:F}).
\end{thm}

\textbf{Sketch of the Proof of Theorem \ref{thm:2}. }  Our approach is based on the application of relative entropy method \cite{Dafermos, Diperna} and more recently \cite{KAN}
, and will be explained in section \ref{sec:proof}. 
The first step of the proof is to introduce a relative entropy which allows us to compare solutions of \eqref{eq:lim} with those of \eqref{eq:feps2}, and then derive an estimate which will enable us to prove that this relative entropy tends to $0$ as $\varepsilon$ vanishes. Deriving such an estimate is the most difficult part of the analysis. More precisely, the difficulties come from: (i) the reaction term $\partial_v(f^\ve N(v))$ which introduces moments of order $4$ that we will need to control, and (ii) the fact that $V^{\varepsilon}$ and $W^{\varepsilon}$ are not a priori bounded in $L^{\infty}(\mathbb{R}^d)$. We will overcome this problem using two entropy inequalities that will be proved in Section \ref{sec:apriori}. 

The rest of the paper is devoted to the proof of Theorem
\ref{thm:2}, where we get {\it a priori} estimates uniformly  with
respect to $\ve>0$ on the solution
$(f^\ve)_{\ve>0}$ constructed in Proposition \ref{prop:exfeps} and
study the behavior of  a relative entropy \cite{Dafermos, Diperna}. 

%% file: apriori.tex
In this section, we prove some a priori estimates of the moments of a
solution to (\ref{eq:feps}) which will be crucial for the proof of
Theorem \ref{thm:2}. For all $i\in\mathbb{N}$ and $z\in\{\mathbf{x},\,v,\,w\}$, we denote by $\mu_i^z$  the moment of order $i$ in $z$ of $f^{\varepsilon}$, defined as
$$
\mu_i^z(t) \,:=\,\int_{\mathbb{R}^{d+2}}|z|^i
f^{\varepsilon}(t,\mathbf{x},v,w)\, \md \mathbf{x} \md v\md w.
$$
whereas  $\mu_i$ is given by
\begin{equation*}
\mu_i(t) \,:=\,\mu_i^\mathbf{x}(t)  \,+\,\mu_i^v(t)  \,+\, \mu_i^w(t).
\end{equation*}

Throughout this sequel, let $T>0$ and $\varepsilon>0$, and suppose
that $f^\ve$ is  a well-defined solution to \eqref{eq:feps} for all
$t\in(0,T]$ obtained in Proposition \ref{prop:exfeps}.  For any $p \in \mathbb{N}^*$, we first establish some a priori estimates on $\mu_{2p}^v$ and $\mu_{2p}^w$.
\begin{lem}[Moment estimates]
\label{lem:M2}
Consider the solution $f^\ve$ to \eqref{eq:feps} given by Proposition
\ref{prop:exfeps} and  $p^*\in \mathbb{N}^*$
$$
\mu^v_{2p^*}(0) + \mu^w_{2p^*}(0) < \infty.
$$ 
Then  there exists $C>0$, only
depending on $p^*$ and $\tau$, such that for all $p\in [1,p^*]$,
\begin{eqnarray}
\label{eq:inM2}
\frac{1}{2p}\,\frac{\md }{\md t}\left[\mu_{2p}^v + \mu_{2p}^w\right](t) \,+\, \mu_{2p+2}^v(t) \,+\,\frac{1}{\ve}\,\mathcal{D}_{p}(t)
\,\leq\, C\,\left(\left[\mu_{2p}^v + \mu_{2p}^w\right](t) \,+\, 1\right),
\end{eqnarray}
where $\mD_p(t)$ is non-negative and  defined as
$$
\mathcal{D}_{p}(t)\, := \,\int_{\mathbb{R}^{d+2}}
f^{\varepsilon}\,v^{2p-1} \left(v-V^{\varepsilon}\right) \,
\rho_0^\ve\,\md \mathbf{x} \md v \md w.
$$
\end{lem}
\begin{Proof}
Consider $f^{\varepsilon}$ a well-defined solution to \eqref{eq:feps}
and $p\in\mathbb{N}^*$, 
we compute the time evolution of moments in $|v|^{2p}$ and $|w|^{2p}$,
hence we have
$$
\dfrac{1}{2p}\,\frac{\md}{\md t}\left[\mu_{2p}^v + \mu_{2p}^w\right](t) \,=\,   I_1 \,+\, I_2 \,+\, I_3 \,+\, I_4 \,+\, I_5,
$$
where
$$
\left\{
\begin{array}{lll}
I_1 &:=&\displaystyle+\int_{\mathbb{R}^{d+2}} v^{2p-1}\,N(v)\, f^{\varepsilon}\,
\md \mathbf{x} \md v \md w, \\ \, \\ 
I_2&:=&\ds-\frac{1}{\varepsilon} \int_{\mathbb{R}^{d+2}}
        f^{\varepsilon}\,\rho^\ve\, v^{2p-1}\,\left(v-V^{\varepsilon}\right)\, \md \mathbf{x} \md v\md w,
\\ \, \\
I_3&:=&\ds-\int_{\mathbb{R}^{d+2}}
  v^{2p-1}\,\K_{\Psi}[f^{\varepsilon}]\,f^{\varepsilon}\, \md \mathbf{x}\md v\md w,
  \\ 
\, \\
I_4 &:=& \displaystyle-\int_{\R^{d+2}}
  v^{2p-1}\,w\,f^{\varepsilon}\, \md \mathbf{x} \md v \md w,
\\ \, \\
I_5&:=&\displaystyle+\int_{\mathbb{R}^{d+2}}
  w^{2p-1}\,A(v,w)\,f^{\varepsilon}\, \md \mathbf{x} \md v \md w.
\end{array}\right.
$$
First of all we treat  the term $I_1$, since $v^{2p-1}\,N(v) \,=\, |v|^{2p}\,-\,|v|^{2p+2}$, we get
$$
I_1\,=\, \displaystyle\int_{\mathbb{R}^{d+2}}
\left(|v|^{2p}-|v|^{2p+2}\right)\,f^{\varepsilon}  \, \md \mathbf{x} \md v \md
w \,=\, \mu_{2p}^v(t) - \mu_{2p+2}^v(t) . 
$$
Furthermore, we estimate the second term $I_2$ and show that it is
non-positive. Indeed, from the definition of $V^\ve$, we
observe that
$$
\int  \rho^\ve \,\left(V^\ve\right)^{2p-1}\,\left(\int  f^{\ve} \,\left(v-V^\ve\right)\,\md v\, \md w \right)\, \md \mathbf{x}=0.
$$ 
Therefore, since $p\geq 1$, the function $v\mapsto
v^{2p-1}$ is non-decreasing,  which yields 
\begin{eqnarray*}
\ve \, I_2 &=& -\int f^{\ve} \,\left[ \left( v^{2p-1} - (V^\ve)^{2p-1}
               \right)  \,+\, (V^\ve)^{2p-1} \right] 
               (v-V^\ve)\,\rho^\ve\,\md v \md w  \md \mathbf{x} 
\\
&=& -\int f^{\ve} \,\left( v^{2p-1} - (V^\ve)^{2p-1}
               \right)\, (v-V^\ve)\, \rho^\ve\, \md v \md w\md \mathbf{x}
       \leq 0.
\end{eqnarray*}
As a consequence, we obtain that $\ve \, I_2 = -\mathcal{D}_p \leq 0$.

Then, we deal with $I_3$ and prove that it contributes to the
dissipation of moments. Indeed, by symmetry of $\Psi$ we may
reformulate $I_3$, using the shorthand notation
$f^{\varepsilon\,\prime}=f^\ve(\mathbf{x}',v',w')$ and
$f^{\varepsilon}=f^\ve(\mathbf{x},v,w)$, as 
\begin{eqnarray*}
 I_3 & =&  -\dfrac{1}{2}\displaystyle\iint \Psi(\mathbf{x}-\mathbf{x}')\,v^{2p-1}\left(v-v'\right)\, f^{\varepsilon\,\prime}(t)\,f^{\varepsilon}(t)\, \md \mathbf{x}\md v \md w\md \mathbf{x}' \md v' \md w' \\
 && - \dfrac{1}{2}\displaystyle\iint \Psi(\mathbf{x}-\mathbf{x}')\,v'^{2p-1}\left(v'-v\right)\, f^{\varepsilon\,\prime}(t)\,f^{\varepsilon}(t)\, \md \mathbf{x}\md v \md w \md \mathbf{x}' \md v' \md w' \\
 & \leq& 0.
\end{eqnarray*}
Finally, we can easily compute $I_4$, which yields  
\begin{eqnarray*}
I_4 & \leq & \int f^\ve \,\left( \dfrac{2p-1}{2p}|v|^{2p} + \dfrac{1}{2p}|w|^{2p} \right)\, \md \mathbf{x} \md v\md w   \\
& \leq& \dfrac{2p-1}{2p}\,\mu_{2p}(t) ,
\end{eqnarray*}
whereas the last term $I_5$ gives
\begin{eqnarray*}
I_5 & =& \tau\, \int \left[ w^{2p-1}\,v \,+\, a\,w^{2p-1} \,-\,
         b\,|w|^{2p}  \right]\,f^{\varepsilon} \, \md \mathbf{x} \md v \md
         w   
\\
& \leq &\tau\,\int \left[\dfrac{1}{2p}\,|v|^{2p}\,+\,\dfrac{2\,p-1}{2\,p}\,|w|^{2p} \,+\, \dfrac{2\,p-1}{2\,p}\,|w|^{2p} \,+\, \dfrac{|a|^{2p}}{2p}\right]\,f^{\varepsilon} \, \md \mathbf{x} \md v \md w\,   \\
& \leq &\tau\,\dfrac{2p-1}{p}\,\mu_{2p}(t) \,+\, \dfrac{\tau |a|^{2p}}{2p}.
\end{eqnarray*}
Therefore, gathering the previous results,  we get the entropy
inequality \eqref{eq:inM2} with $C>0$, only depending on $p$, $\tau$ and $a$.
\end{Proof}

This Lemma will be helpful to pass to the limit $\ve\rightarrow 0$ in
(\ref{eq:feps}). The first consequence is some moment  estimates in
$\mathbf{u}=(v,w)$ and also $\mathbf{x}$.
\begin{cor}
\label{cor:M2}
Under the assumptions of Lemma \ref{lem:M2} with $p^*=2$, we choose an
initial datum $f^\ve_0$ such that \eqref{H2-f0:M4}
is satisfied. Then, there exists a constant $C_T> 0$, which does not
depend on $\ve>0$,
such that for any $k\in[0, 4]$,
\begin{equation} 
\label{estim:Moments}
\left\{
\begin{array}{l}
\ds \mu_{k}(t)  \,\leq\, C_T, \quad t\in [0,T],
\\ \, \\
\ds\int_0^T \mu_{k+2}^v(t) \, \md t \, \leq \, C_T.
\end{array}\right.
\end{equation}
\end{cor}
\begin{Proof}
First we observe that 
\begin{equation}
\label{r:1}
\mu_0(t) \,=\, 3\, \int_{\mathbb{R}^{d+2}} f^\ve(t) \,\md
v \md w \md \mathbf{x} \,=\, 3\, \|\rho_0^\ve\|_{L^1},  
\end{equation}
which gives the result with $k=0$. 

Then for $k=4$, we apply Lemma \ref{lem:M2} with $p=2$ and integrate with respect to $t\in
[0,T]$ and by the Gr\"onwall's lemma, it yields the estimates on the
second order moment in $\mathbf{u}=(v,w)$, there exists a constant $C_T>0$ such
that for any $t\in [0,T]$
$$
 \mu_{4}^v(t)\, +\, \mu_{4}^w(t) \leq C_T.
$$
Furthermore, since $\rho^\ve$ does not depend on time, we also have
for all $t\in [0,T]$,
$$
\mu_{4}^\mathbf{x}(t) \,=\, \int_{\mathbb{R}^d} |\mathbf{x}|^4
\,\rho_0^\ve \,\md \mathbf{x} \,< \,\infty,   
$$ 
hence from hypothesis \eqref{H2-f0:M4} which gives us the uniform control of the moment of order $4$ in $\mathbf{x}$ of $f_0^\ve$, there
exists a constant $C_T>0$, independent of $\ve>0$, such that for all $t\in [0,T]$,
\begin{equation}
\label{r:2}
\mu_{4}(t) \,\leq\, C_T.
\end{equation}
On the other hand, from the latter result and the dissipative terms
obtained in Lemma \ref{lem:M2}, there exists a constant $C_T>0$ such
that 
\begin{equation}
\label{r:3}
\int_0^T \mu_{6}^v(t) \, \md t \,\leq \, C_T.
\end{equation}
Interpolating (\ref{r:1}) and
(\ref{r:2})-(\ref{r:3}), it yields the result with $0\leq k\leq 4$.
\end{Proof}

Another consequence of Lemma \ref{lem:M2} is the control of the
dissipation $\mathcal{D}_{1}(.)$, which is a crucial step to characterize
the limit of the sequence $(f^\ve)_{\ve>0}$ when $\ve \rightarrow 0$. 
\begin{cor}\label{corD1}
Under the assumptions of Lemma \ref{lem:M2} with $p^*=1$, we choose an
initial datum $f^\ve_0$ such that \eqref{H2-f0:M4}
is satisfied. Then, there exists a constant $C_T> 0$,
such that
\begin{equation}
\int_0^T\int_{\R^{d+2}} f^{\varepsilon}(t)\,|v-V^{\varepsilon}(t,\mathbf{x})|^2
\, \rho^\ve_0(\mathbf{x})\,\md \mathbf{x}\, \md v\, \md w\, \md t \,\leq\, C_T\,\varepsilon.
\label{estimate1}
\end{equation}
\end{cor}
\begin{Proof}
We first apply Corollary \ref{cor:M2} to obtain a uniform estimate on
the moments $\mu_2(.)$,
$$
\mu_{2}(t)  \,\leq\, C_T, \quad t\in [0,T].
$$
Then, integrating the entropy inequality (\ref{eq:inM2}) between $0$ and $T$ and using the positivity of $\mathcal{D}_1(.)$, we find
\begin{eqnarray*}
&&\frac{1}{\varepsilon}\int_0^T\int_{\R^{d+2}}
  f^{\varepsilon}(t,\mathbf{x},v,w)\,|v-V^{\varepsilon}(t,\mathbf{x})|^2 \, \rho_0^\ve(\mathbf{x})\md \mathbf{x} \,\md v\, \md w \,\md t \\
&& \leq\,    \mu_2(0)  \,+\, 2\,C \,\int_0^T \left(1\,+\,\mu_2(t) \right)\md t \, \leq\, C_T,
\end{eqnarray*}
from which we easily deduce \eqref{estimate1}.
\end{Proof}

This last result is not enough to justify the asymptotic limit. Hopefully, it can be improved by removing the weight  $\rho^\ve_0$
in the previous estimate.
\begin{lem}\label{corD1sansrho}
Consider the solution $f^\ve$ to \eqref{eq:feps} given by Proposition
\ref{prop:exfeps}, where the initial datum $f^\ve_0$ satisfies \eqref{H2-f0:M4}.
Then, there exists a constant $C_T>0$, independent of $\ve>0$, such that
\begin{equation}
\int_0^T\int_{\R^{d+2}}
f^{\varepsilon}(t,\mathbf{x},v,w)|v-V^{\varepsilon}(t,\mathbf{x})|^2\, \md \mathbf{x}\, \md v\, \md w\, \md t \leq C_T\, \varepsilon^{{2}/{(d+6)}}.
\label{estimate1sansrho}
\end{equation}
\end{lem}
\begin{Proof}
Let us fix $\ve>0$ and $T>0$. First, we set
$$
I^\ve \,=\,  \int_0^T\int f^{\varepsilon} |v-V^{\ve}|^2\, \md w\, \md v\, \md \mathbf{x}\, \md t
$$
and notice that according to the definition of $V^\ve$,
\begin{eqnarray*}
I^\ve & =& \displaystyle\int_0^T\int f^\ve \left( |v|^2 + |V^\ve|^2 -
           2vV^\ve \right) \md w\, \md v\, \md \mathbf{x}\, \md t  
 \\
&  \leq& \displaystyle\int_0^T\int f^{\varepsilon} |v|^2 \md w\, \md v\, \md \mathbf{x}\, \md t  \,<\, +\infty.
\end{eqnarray*}
so it gives that $f^{\varepsilon} |v-V^{\ve}|^2\in
L^1([0,T]\times\mathbb{R}^{d+2})$. Our strategy to prove \eqref{estimate1sansrho} is to divide $\R^{d}$ into several subsets where $f^{\varepsilon}
|v-V^{\varepsilon}|^2$ is easier to control. 

We consider any $\eta>0$ and  define the set $\mathcal{A}_\ve$
$$
\mathcal{A}_\ve:=\left\{ \mathbf{x}\in\mathbb{R}^d \quad|\quad
  \rho_0^\ve(\mathbf{x})=0 \right\},
$$
and  $\mathcal{B}^{\eta}_\ve$ given by 
$$
\mathcal{B}^{\eta}_\ve:=\left\{ \mathbf{x}\in\mathbb{R}^d \quad|\quad
  \rho_0^\ve(\mathbf{x})>\eta \right\},
$$
whereas  $\mathcal{C}^{\eta}_\ve=\mathbb{R}^d \backslash \left(A _\ve\cup
  B^{\eta}_\ve  \right)$, that is, 
\begin{eqnarray*}
\mathcal{C}^{\eta}_\ve&:=& \left\{ \mathbf{x}\in\mathbb{R}^d \quad|\quad
    0<\rho_0^\ve(\mathbf{x})\leq\eta \right\}.
\end{eqnarray*}
Thus, we have $ I^\ve \,=\,  I_1^\ve + I_2^\ve + I_3^\ve$,  where
$$
\left\{
\begin{array}{l}
I_1^\ve:=\displaystyle\int_0^T\int_{\mathcal{A}_\ve}\int f^{\ve} |v-V^{\ve}|^2
  \md w\, \md v\, \md \mathbf{x}\, \md t, \\ \,
\\ 
I_2^\ve:=\displaystyle\int_0^T\int_{\mathcal{B}^{\eta}_\ve}\int f^{\ve}
  |v-V^{\ve}|^2 \md w\, \md v\, \md \mathbf{x}\, \md t,
\\ \, \\
I_3^\ve:=\displaystyle\int_0^T\int_{\mathcal{C}^{\eta}_\ve}\int f^{\ve} |v-V^{\ve}|^2  \md w\, \md v\, \md \mathbf{x}\, \md t.
\end{array}\right.
$$
On the one hand, since $f^{\ve}\geq 0$, we directly have that $f^\ve =
0$ when $\rho_0^\ve=0$, that is when $\mathbf{x} \in \mathcal{A}_\ve$, thus we have 
\begin{equation}\label{I1}
I_1^\ve = \displaystyle\int_0^T\int_{\mathcal{A}_\ve}\int f^{\ve} |v-V^{\ve}|^2 \md \mathbf{x} \md v \md w \md t = 0. 
\end{equation}
On the other hand,  for $\mathbf{x}\in \mathcal B_\ve^\eta$, we know that
$\rho^\ve_0(\mathbf{x})> \eta$, therefore it yields that 
\begin{eqnarray*}
I_2^\ve &\leq& \displaystyle\int_0^T\int_{\mathcal{B}^{\eta}_\ve}\int f^{\ve}\,
  |v-V^{\ve}|^2 \,\dfrac{\rho_0^\ve}{\eta}\, \md \mathbf{x} \,\md v\, \md w\, \md t
\\
 &\leq&
  \dfrac{1}{\eta}\displaystyle\int_0^T\int
  f^{\ve} \,|v-V^{\ve}|^2 \,\rho_0^\ve\,\md \mathbf{x} \,\md v \,\md w \,\md t. 
\end{eqnarray*}
Hence, by application of Corollary \ref{corD1},  we get the following estimate: 
\begin{equation}\label{I2}
I_2^\ve = \mathcal{O}\left(\dfrac{\ve}{\eta}\right).
\end{equation}

It remains to control the last term $I_3^\ve$. To this aim we bound  it by
the sum of three terms : for any $R>0$, we have
$$
I_3^\ve  \,\leq\,\int_0^T\int_{\mathcal{C}^{\eta}_\ve}\int f^{\ve} |v|^2 \md v\, \md \mathbf{x}\, \md w\,\md t    =:  I_{3,1}^\ve \,+\; I_{3,2}^\ve \,+\, I_{3,3}^\ve,
$$
where 
$$
\left\{
\begin{array}{l}
I_{3,1}^\ve \,:=\,\displaystyle\int_0^T\int_{\mathcal{C}^{\eta}_\ve}\int_{\{ |v|> R
  \}} f^{\ve} |v|^2 \md v \,\md \mathbf{x}\, \md w\,\md t, \\ \, \\
 I_{3,2}^\ve \,:=\,\displaystyle\int_0^T\int_{\mathcal{C}^{\eta}_\ve \cap
  B^c(0,R)}\int_{\{ |v|\leq R \}} f^{\ve} |v|^2 \md v \,\md \mathbf{x}
  \,\md w\,\md t,
\\ \, \\
I_{3,3}^\ve \,:=\,\displaystyle\int_0^T\int_{\mathcal{C}^{\eta}_\ve \cap
  B(0,R)}\int_{\{ |v|\leq R \}} f^{\ve} |v|^2 \md v \,\md \mathbf{x}
  \,\md w\,\md t.
\end{array}\right.
$$
For $k>2$, the estimates on moments in velocity gives that
\begin{equation}
\label{I31}
I_{3,1}^\ve  \,\leq\, \frac{1}{R^{k-2}}\int_0^T\int f^\ve |v|^k \,\md v
\,\md \mathbf{x}\, \md w \,\md t \,=\, \dfrac{1}{R^{k-2}}\int_0^T \mu_k^v(t) \,\md t ,
\end{equation}
where the last term is uniformly bounded according to Corollary \ref{cor:M2} for
$k\in [0,6]$.  

Furthermore,  we estimate $I_{3,2}^\ve$ using a similar
argument as for $I_{3,1}^\ve $ but now using moments in space, that is, for $p\in[2,4]$, 
\begin{equation}
\label{I321}
I_{3,2}^\ve  \leq \int_0^T\int \dfrac{|\mathbf{x}|^p}{R^p}\, f^{\ve}(t)\, R^2 \,\md w \,\md v \,\md \mathbf{x} \,\md t \,\leq \,\dfrac{1}{R^{p-2}} \, \int_0^T \mu_p^\mathbf{x}(t) \,\md t ,
\end{equation}
where $\mu_p^\mathbf{x}(.)$ is uniformly bounded according  to Corollary \ref{cor:M2} for
$p\in [0,4]$. Finally, the last term $I_{3,3}^\ve$  can be  computed as
\begin{equation}\label{I322}
I_{3,3}^\ve  \,\leq\,R^2 \;\int_0^T\int_{\mathcal{C}^{\eta}_\ve \cap B(0,R)} \rho_0^{\ve} \, \md \mathbf{x} \,\md t \,\leq\; \tilde{C}\,T\,R^{d+2}\,\eta,
\end{equation}
where $\tilde{C}$ is the positive constant such that $|B(0,R)| =
\tilde{C} R^d$.  By summing \eqref{I31}, \eqref{I321} and \eqref{I322}, we can conclude that there exists a constant $C=C(T)$ independent on $R$, $\eta$ and $\varepsilon$ such that
$$
I_3 ^\ve  \,\leq\, C\,\left(\, R^{d+2}\,\eta \,+\, \dfrac{1}{R^{p-2}} \,+\, \dfrac{1}{R^{k-2}} \,\right).
$$
For simplicity we  choose $p=k=4$ and optimize the value of $R$, which
leads to
\begin{equation}
\label{I3}
I_3 ^\ve  \,\leq\, C_T \,\eta^{{2}/{(d+4)}}.
\end{equation}

Finally by summing \eqref{I1}, \eqref{I2} and \eqref{I3}, we  get that there exists a positive constant $C_T>0$ independent on $\eta$ and $\ve$, such that,
$$
I^\ve \,\leq\, C_T\, \left( \eta^{{2}/{(d+4)}} \,+\, \dfrac{\ve}{\eta}  \right).
$$
Hence, we can choose  $\eta = \ve^{{(d+4)}/{(d+6)}}$ and
$$
\displaystyle\int_0^T\int f^{\varepsilon} |v-V^{\ve}|^2 \,\md \mathbf{x} \,\md v\,
\md w \, \md t \,\leq\, C_T\, \ve^{{2}/{(d+6)}}.
$$
\end{Proof}

%% file: entropy.tex
Following ideas from \cite{FIG,KAN}, our proof of Theorem \ref{thm:2} relies on a relative entropy argument  to estimate the distance between a solution of (\ref{eq:feps2}) and a solution of (\ref{eq:lim}) on $[0,T]$ for some finite $T>0$ as $\varepsilon\rightarrow0$, with well-prepared initial conditions. Throughout this section, by a solution $(\rho_0,\rho_0 V,\rho_0 W)$ of (\ref{eq:lim}), we refer to the solution constructed in Proposition \ref{prop:exhydro}. We start this section with the definition of the relative entropy that we will be using. Then, we present an equality satisfied by the relative entropy, which will be useful to estimate it. This estimate will finally be the key argument to conclude the proof.

\subsection{Definition of relative entropy}

We want to use a relative entropy argument which enables us to compare solutions of (\ref{eq:lim}) with the solutions of (\ref{eq:feps}) seen on its hydrodynamic form (\ref{eq:feps2}). We first introduce the notion of entropy.

\begin{defi}[Entropy]\label{defi:entr}
For all functions $V$ and $W:\mathbb{R}^{d}\longrightarrow\mathbb{R}$, and for any non-negative function $\rho:\mathbb{R}^d\rightarrow\mathbb{R}$, we define for $\mathcal{Z}=(\rho,\rho V, \rho W)$, the  entropy $\eta(\mathcal Z)$ by 
\begin{equation}
\label{entropy}
\eta(\mathcal Z)\,:=\,\rho\, \frac{|V|^2 \,+\, |W|^2}{2}.
\end{equation}
\end{defi}
Note that if we define $P=\rho V$ and $Q=\rho W$, then we have $$
\eta(\mathcal Z)\,=\,\frac{P^2\,+\,Q^2}{2\,\rho}.
$$
As a consequence, the differential of $\eta$ with respect to its variables $\mathcal Z$ is given by
\begin{equation}
\label{Deta}
D\eta(\mathcal Z) \,=\,  \begin{pmatrix} D_\rho \eta \\  D_P \eta \\  D_Q\eta \end{pmatrix}=  \begin{pmatrix}\ds -\,\frac{|V|^2 + |W|^2}{2} \\  V \\  W \end{pmatrix}.
\end{equation}
%

Using the definition of entropy \eqref{entropy}, we can now introduce the notion of relative entropy.

\begin{defi}[Relative Entropy]\label{defi:entrrel}
For all functions $V_1$, $W_1$, $V_2$ and $W_2:\mathbb{R}^{d}\to\mathbb{R}$ and for all non-negative functions $\rho_1$ and $\rho_2:\mathbb{R}^d\to \mathbb{R}$, we define for $\mathcal Z_i=(\rho_i,\rho_iV_i,\rho_i W_i)$, with  $i=1,2$, the relative entropy 
$$
\eta(\mathcal Z_1| \mathcal Z_2)  \,:=\,  \eta(\mathcal Z_1) \,-\, \eta(\mathcal Z_2)   \,-\, D\eta(\mathcal Z_2)\cdot \left( \mathcal Z_1 - \mathcal Z_2\right)  
$$
which gives us after computation
\begin{equation}\label{eq:entrel}
\eta(\mathcal Z_1| \mathcal Z_2) \,=\, \rho_1 \,\displaystyle\frac{|V_2-V_1|^2 \,+\, |W_2-W_1|^2}{2}.
\end{equation}
\end{defi}
This relative entropy will be useful to "compare" the weak solution $(\rho_0^{\varepsilon},\rho_0^{\varepsilon}V^{\varepsilon},\rho_0^{\varepsilon}W^{\varepsilon})$ to (\ref{eq:feps2}) with the classical solution $(\rho_0,\rho_0 V,\rho_0 W)$ to (\ref{eq:lim}).

\subsection{Relative entropy equality}

In this subsection, we prove an equality satisfied by the relative entropy defined by (\ref{eq:entrel}) between an arbitrary  smooth function and a solution of the hydrodynamic equations (\ref{eq:lim}). The purpose of this result is to split the relative entropy dissipation into one part due to the macroscopic solution  (\ref{eq:lim}) and another part which estimates the difference between the two solutions. 
\begin{lem}
Under the assumption that $\Psi$ is non-negative, symmetric and satisfies
(\ref{H1}), we consider $(V,W)$ the solution to the hydrodynamic equation
(\ref{eq:nlFHN}) given by Proposition \ref{prop:exhydro}. Then, for any  $\widetilde{\mathcal
  Z}\,=\,(\widetilde{\rho}, \,\widetilde{\rho}\,\widetilde{V},
\,\widetilde{\rho}\,\widetilde{W})$ such that $\widetilde{\rho}$ is
non-negative, $\widetilde{\rho}\in L^1\cap
L^{\infty}(\mathbb{R}^d)$,  whereas $\widetilde{V}$ and
$\widetilde{W}$ are both differentiable in time and such that for any
$t\in[0,T]$,  
$$
\widetilde{\rho} \left(|\widetilde{V}(t)|^4 + |\widetilde{W}(t)|^4
\right) \in L^1(\mathbb{R}^d),
$$
the following equality holds:
\begin{eqnarray}
\label{eq:entropy}
\ds\frac{\md}{\md t}\ds\int \eta(\widetilde{\mathcal{Z}}|\mathcal{Z})\,\md \mathbf{x}\, 
 &=&\displaystyle\frac{\md}{\md t}\displaystyle\int
     \eta(\widetilde{\mathcal{Z}})\md \mathbf{x} 
\\
\nonumber 
&&- \displaystyle\int D\eta(\mathcal{Z})\left [ \partial_t \widetilde{\mathcal{Z}} - \mF(\widetilde{\mathcal{Z}}) \right ] \md \mathbf{x} 
\,+\,\, \mR(\widetilde\mZ|\mZ) \,+\, \mS(\widetilde\mZ), 
\end{eqnarray}
with $\mR \,=\, \mR_l +\mR_{nl}$ and $\mS\,=\,\mS_l+\mS_{nl}$, where $\mR_l$ and $\mS_l$ contain local terms 
\begin{equation}
\label{def:RS}
\left\{
\begin{array}{lll}
\mR_l(\widetilde\mZ|\mZ) &:=& \ds\int \widetilde{\rho}\left(V-\widetilde{V}\right)\,\left(N(V)-W \,-\, N(\widetilde{V})+\widetilde{W} \right) \,\md \mathbf{x}\\ \, \\
&\,&+\,\ds\int \widetilde{\rho}\left(W-\widetilde{W}\right)\,\left(A(V,W)-A(\widetilde{V},\widetilde{W})\right)\,\md \mathbf{x},
\\ \, \\
\mS_l(\widetilde\mZ) &:=& -\displaystyle\int \widetilde{\rho}\,\left[\widetilde{V}\,N(\widetilde{V})\;-\, \widetilde{V}\,\widetilde{W}\, \,+\,\widetilde{W}\,A(\widetilde{V},\widetilde{W}) \right]\,\md \mathbf{x} .
\end{array}\right.
\end{equation}
whereas $\mR_{nl}(\widetilde\mZ|\mZ)$ and $\mS_{nl}(\widetilde\mZ)$ gather nonlocal terms 
\begin{equation}
\label{def:RS2}
\left\{
\begin{array}{lll}
\mR_{nl}(\widetilde\mZ|\mZ) &:=& \ds\int \tilde \rho \,\left[(V-\widetilde{V})
  \,\left(\mathcal{L}_{\rho_0}(V)
    -\mathcal{L}_{\widetilde\rho}(\widetilde V)\right)\right] \md \mathbf{x},
\\ \, \\
\mS_{nl}(\widetilde\mZ) &:=& \ds-\frac{1}{2}\iint \Psi(\mathbf{x}-\mathbf{y})\,\widetilde{\rho}(\mathbf{x})\,\widetilde{\rho}(\mathbf{y}) \left|\widetilde{V}(t,\mathbf{x}) \,-\, \widetilde{V}(t,\mathbf{y})\right|^2 \,\md \mathbf{x}\, \md \mathbf{y} 
\end{array}\right.
\end{equation}
\label{lem:gen}
\end{lem}
\begin{Proof}
First of all, it is worth noticing that for all $t\in[0;T]$, since
$V(t)$ and $W(t)\in L^{\infty}(\mathbb{R}^d)$ and since $\rho_0\in
L^1(\mathbb{R}^d)$, then for all $i\in\mathbb{N}$,
$\rho_0\left(|V(t)|^i+|W(t)|^i\right) \in L^1(\mathbb{R}^d)$. From the
definition (\ref{defi:entrrel}) of $\eta(\widetilde{\mathcal{Z}}|\mathcal{Z})$, we have
$$\begin{array}{l l l}
 \displaystyle\frac{\md}{\md t} \displaystyle\int \eta(\widetilde{\mathcal{Z}}|\mathcal{Z})\, \md \mathbf{x} & = & \displaystyle\int \left[ \partial_t \eta(\widetilde{\mathcal{Z}}) - \partial_t \eta(\mathcal{Z}) - \partial_t D\eta(\mathcal{Z})\cdot\left( \widetilde{\mathcal{Z}}-\mathcal{Z} \right)  - D\eta(\mathcal{Z})\cdot\partial_t\left( \widetilde{\mathcal{Z}}-\mathcal{Z} \right) \right]\md \mathbf{x}\\
 & = & I_1 + I_2,
\end{array}$$
with 
$$
\left\{
\begin{array}{l l l}
I_1 & := & \displaystyle\int \partial_t \eta(\widetilde{\mathcal{Z}})\, -\, D\eta(\mathcal{Z})\cdot\left[\partial_t \widetilde{\mathcal{Z}} \,-\, \mF(\widetilde\mZ)\right]\,\md \mathbf{x},
\\ \, \\
I_2 &:=&
- \displaystyle\int \left[\partial_t D\eta(\mathcal{Z})\cdot\left( \widetilde{\mathcal{Z}}-\mathcal{Z} \right)  \,+\,  D\eta(\mathcal{Z})\cdot \mF(\widetilde{\mZ})\right] \,\md \mathbf{x},
 \end{array}\right.
$$
where $\mF$ is defined in (\ref{eq:lim})-(\ref{def:F}).

On the one hand, the  term $I_1$ corresponds to the variation of
entropy which simply gives 
\begin{equation}
\label{L:I1}
I_1 \,=\, \frac{\md }{\md t }\int \eta(\widetilde{\mathcal{Z}}) \,\md \mathbf{x} \,-\,\int D\eta(\mathcal{Z})\cdot\left[\partial_t \widetilde{\mathcal{Z}} \,-\, \mF(\widetilde\mZ)\right]\,\md \mathbf{x}.
\end{equation}

On the other hand, we decompose $I_2$ as   $I_2=I_{21}+I_{22}$ with 
$$
\left\{
\begin{array}{l}
I_{21} \,:=\,\ds\int \partial_t D\eta(\mathcal{Z})\cdot\left( \widetilde{\mathcal{Z}}-\mathcal{Z} \right) \,\md \mathbf{x},
\\ \, \\ 
I_{22} \,:=\, - \displaystyle\int D\eta(\mathcal{Z})\,\cdot \,\mF(\widetilde{\mathcal{Z}}) \md \mathbf{x}. 
\end{array}\right.
$$
Using the definition of $D\eta(\mZ)$ in \eqref{Deta} and since $\mathcal{Z}$ is solution to (\ref{eq:lim}), we have
\begin{eqnarray*}
I_{21} & := & - \displaystyle\int\partial_t \,\begin{pmatrix} \ds- \frac{|V|^2 + |W|^2}{2}  \\  V  \\  W  \end{pmatrix} \cdot \begin{pmatrix} \widetilde{\rho} \,-\, \rho_0  \\  \widetilde{\rho}\,\widetilde{V} \,-\, \rho_0 \,V  \\ \widetilde{\rho}\,\widetilde{W} \,-\, \rho_0 \,W  \end{pmatrix}\, \md \mathbf{x}   \\
& = & \displaystyle\int\left(\widetilde{\rho} \,-\, \rho _0 \right)\,\left[ V\,\left(\mathcal{L}_{\rho_0}(V)\,+\,N(V)\,-\,W\,\right) \,+\, W\,A(V,W) \right] \,\md \mathbf{x} \\ 
& & - \, \displaystyle\int\left(\widetilde{\rho}\,\widetilde{V} - \rho_0 V\right)\left(\mathcal{L}_{\rho_0}(V)+N(V)-W\right) \,\md \mathbf{x}  \\
& & - \, \displaystyle\int\left(\widetilde{\rho}\,\widetilde{W} \,-\, \rho_0\, W\right)\,A(V,W)\,\md \mathbf{x},  \\
\end{eqnarray*}
hence it  yields
\begin{equation*}
 I_{21} \,=\, \int\tilde\rho \,\left[\left(V\,-\,\widetilde{V}\right)\,\left[\mathcal{L}_{\rho_0}(V)\,+\,N(V)\,-\,W\right] \,+\,A(V,W)\,\left(W\,-\,\widetilde{W}\right) \right]\,\md \mathbf{x}.
\end{equation*}
Furthermore, from the definition of  $\mF(\widetilde{\mathcal{Z}})$ in \eqref{eq:lim}-\eqref{def:F} and $D\eta(\mathcal{Z})$ in \eqref{Deta}, we obtain 
\begin{equation*}
 I_{22} \,=\, -\int \widetilde{\rho} \,\left( V\,\left[\mathcal{L}_{\tilde\rho}(\widetilde{V})\,+\,N(\widetilde{V})\,-\,\widetilde{W}\right]\,+\,W \, A(\widetilde{V},\widetilde{W})\right)\,\md \mathbf{x}.
\end{equation*}
Then, gathering the latter two equalities and after reordering, we have
\begin{eqnarray*}
I_2 &=& \int \tilde \rho \,\left[(V-\widetilde{V}) \,\left(\mathcal{L}_{\rho_0}(V) -\mathcal{L}_{\widetilde\rho}(\widetilde V)\right) \,-\, \widetilde V\,\mathcal{L}_{\tilde\rho} (\widetilde{V})\right] \md \mathbf{x}\,
\\
& & + \, \int \widetilde{\rho} \,\left[\left(V \,-\, \widetilde{V}\right)\,\left(N(V) \,-\; N(\widetilde{V})\,-\, W \,+\, \widetilde{W}\right) \;-\, \widetilde{V}\,N(\widetilde{V}) \,+\, \widetilde{V}\,\widetilde{W}   \right]\,\md \mathbf{x}
\\
& & + \,\int \widetilde{\rho} \,\left[ \left(W\,-\,\widetilde{W}\right)\,\left(A(V,W)\,-\,A(\widetilde{V},\widetilde{W})\right)\,-\,\widetilde{W}\,A(\widetilde{V},\widetilde{W})\, \right]\,\md \mathbf{x},
\end{eqnarray*}
which can be also written as 
\begin{equation}
\label{L:I2}
I_2 \,=\ \mI_{\psi}\,+\, \mR_l(\widetilde{\mZ}|\mZ) \,+\, \mS_l(\widetilde{\mZ}),
\end{equation}
where $\mR_l$ and $\mS_l$ are given in \eqref{def:RS} whereas the
first term $\mI_\psi$ is given by
$$
\mI_\psi\,:=\,\int \tilde \rho \,\left[(V-\widetilde{V})
  \,\left(\mathcal{L}_{\rho_0}(V)
    -\mathcal{L}_{\widetilde\rho}(\widetilde V)\right) \,-\,
  \widetilde V\,\mathcal{L}_{\tilde\rho} (\widetilde{V})\right] \md \mathbf{x}.
$$
Thus, we set 
$$
\left\{
\begin{array}{l}
\ds\mR_{nl} \,=\, \int \tilde \rho \,\left[(V-\widetilde{V})
  \,\left(\mathcal{L}_{\rho_0}(V)
    -\mathcal{L}_{\widetilde\rho}(\widetilde V)\right)\right] \md \mathbf{x},
\\ \,\\
\ds\mS_{nl}  \,=\, -\int \tilde \rho \,\widetilde V\,\mathcal{L}_{\tilde\rho}(\widetilde{V})\md \mathbf{x}
\end{array}\right.
$$
and a direct computation gives 
$$
\mS_{nl}  \,=\,
-\displaystyle\frac{1}{2}\,\iint\Psi(\mathbf{x}-\mathbf{y})\,\widetilde{\rho}(\mathbf{x})\,\widetilde{\rho}(\mathbf{y})\,\left|
  \widetilde{V}(t,\mathbf{x}) - \widetilde{V}(t,\mathbf{y}) \right|^2 \,\md \mathbf{x}\,\md \mathbf{y}.
$$

We conclude the proof by gathering the last equality together with \eqref{L:I1} and \eqref{L:I2}. 
\end{Proof}

\subsection{Relative entropy estimate}

Now that we have established the relative entropy equality (\ref{eq:entropy}), we apply it with $\widetilde{\rho}=\rho^{\varepsilon}_0$, $\widetilde{V}=V^{\varepsilon}$ and $\widetilde{W}=W^{\varepsilon}$ to estimate the relative entropy between the weak solution $(V^{\varepsilon},W^{\varepsilon})$ and the classical solution $(V,W)$. More precisely, we prove the following result.
\begin{prop}\label{prop:estimentrop}  
Under the assumptions of Theorem \ref{thm:2}, there exists $C_T>0$ such that we have for all $t\in (0,T]$:
\begin{equation}
\label{eq:estimentrop}
\int_{\R^d} \rho_0^{\varepsilon}(\mathbf{x})\left(|V^{\varepsilon}(t,\mathbf{x}) - V(t,\mathbf{x})|^2 +
  |W^{\varepsilon}(t,\mathbf{x}) - W(t,\mathbf{x})|^2\right) \md \mathbf{x}  \,\leq\, C_T \,\ve^{1/(d+6)},
\end{equation}
where $(V,W)$ is  the solution to (\ref{eq:nlFHN}) and $(\rho_0^\ve, \rho_0^\ve V^\ve, \rho_0^\ve W^\ve)$
are the macroscopic quantities computed from $f^{\varepsilon}$ the solution to (\ref{eq:feps}) on $[0,T]$.
\end{prop}
\begin{Proof} Consider $f^\ve$ a solution to (\ref{eq:feps}) on
  $[0,T]$ given in Proposition \ref{prop:exfeps}. From Corollary \ref{cor:M2}, we get that for any $t\in
  [0,T]$, the moment $\mu_4(t)$  is uniformly bounded with respect to
$\ve>0$.  Therefore applying the H\"older inequality, we obtain that for all $\mathbf{x}\in\R^d$ such that $\rho^\ve_0(\mathbf{x})>0$ and for all $t\in[0;T]$,
\begin{eqnarray}
\nonumber
\rho_0^{\varepsilon}(\mathbf{x}) |V^{\varepsilon}(t,\mathbf{x})|^{4} &=& \displaystyle\frac{1}{|\rho_0^\ve(\mathbf{x})|^3}\,\left( \int vf^{\varepsilon}(t,\mathbf{x},v,w) \md v \, \md w \right)^{4}, \\
&\leq& \int |v|^{4}\,f^\ve(t,\mathbf{x},v,w)  \,\md v\, \md w.
\label{estimatesVeps24}
\end{eqnarray}
 Note that the last inequality remains true when
 $\rho^\ve_0(\mathbf{x})=0$ and the same argument applies when we replace
 $V^\ve$ by  $W^{\varepsilon}$. 

Consequently, since $\rho^\ve \in L^1(\mathbb{R}^d)$, we get for any
$0\leq p \leq 4$ and  $t\in[0,T]$, 
$$
\rho_0^{\varepsilon} \left( |V^{\varepsilon}(t)|^p\,+\,
  |W^{\varepsilon}(t)|^p\right) \,\in
\;L^1(\mathbb{R}^d).
$$
Thus, we can compute the time evolution of the entropy $\eta(\mZ^\ve)$
where $\mZ^\ve=(\rho_0^\ve, \rho_0^\ve V^\ve, \rho_0^\ve W^\ve)$
corresponds to the moments of $f^\ve$ with respect to $(1,v,w)$, that
is $\mZ^\ve$ is solution to \eqref{eq:feps2}. It
yields that
\begin{equation}\label{etaZeps}
\frac{\md}{\md t}\displaystyle\int \eta(\mZ^\ve(t))\,\md \mathbf{x}
\,+\, \mS(\mZ^\ve(t)) \,=\, \int\,V^\ve(t)\,\mathcal{E}(f^\ve(t)) \, \md \mathbf{x}, 
\end{equation}
where $\mS(\mZ^\ve)$ is given by \eqref{def:RS}-\eqref{def:RS2} when
$\widetilde\mZ=\mZ^\ve$ and the error term $\mathcal{E}(f^\ve)$ is defined in \eqref{errorR}.

Then we  consider  $(V,W)$ the solution to (\ref{eq:nlFHN}) given in Proposition
  \ref{prop:exhydro}, hence we have for $\mZ=(\rho_0, \rho_0\,V, \rho_0\, W)$,
\begin{equation}
\label{DetaZeps}
\int D\eta( \mZ(t) )\, [ \partial_t\mZ^\ve(t) - \mF( \mZ^\ve(t) ) ]\,\md \mathbf{x}\, \, = \, \int V(t)\,\mathcal{E}(f^\ve(t))\md \mathbf{x}\,
\end{equation}

Therefore, applying Lemma \ref{lem:gen} with
$\widetilde \mZ=(\rho_0^\ve,\rho_0^\ve\,V^\ve,\rho_0^\ve\,W^\ve)$ and using \eqref{etaZeps} and \eqref{DetaZeps}, we simply get the following equality
\begin{eqnarray}
\label{dd}
\frac{d}{dt}\int \eta({\mathcal{Z}}^\ve(t)|\mathcal{Z}(t))\,\md \mathbf{x}
  \,=\,  
\int \left(V^\ve-V\right)\,\mathcal{E}(f^\ve)(t,\mathbf{x})\md \mathbf{x}
  \,+\, \mR(\mZ^\ve|\mZ), 
\end{eqnarray}
where $\mR=\mR_l + \mR_{nl}$ is given in
(\ref{def:RS})-(\ref{def:RS2}). On the one hand,   we  estimate the
term $\mR_l $ by
\begin{eqnarray*}
\left|\mR_l(\mZ^\ve|\mZ)\right| & \leq& \left|\int \rho_0^\ve\,\left( V-V^\ve \right)\,\left(
  N(V)-N(V^\ve)\, -\, (W-W^\ve) \right)\,\md \mathbf{x} \right| 
\\
& & + \,\left|\int \rho_0^\ve\left( W-W^\ve \right)\,\left( A(V,W)
  \,-\,A(V^\ve,W^\ve) \right) \,\md \mathbf{x} \right|  
\\
& \leq& \int \left( |V-V^\ve|^2\, + \, \dfrac{1}{2}\,(
  |V-V^\ve|^2+|W-W^\ve|^2 ) \right)\, \rho_0^\ve\,\md \mathbf{x}  
\\
& & + \,\dfrac{\tau}{2}\,\int\left( |V-V^\ve|^2+|W-W^\ve|^2 \right) \,\rho_0^\ve\,\md \mathbf{x}   \\
&\leq& \dfrac{3+\tau}{2} \,\int\eta(\mZ^\ve|\mZ)\md \mathbf{x}.
\end{eqnarray*}
On the other hand, we estimate the second term $\mR_{nl}(\mZ^\ve|\mZ)$ as
\begin{eqnarray*}
|\mR_{nl}(\mZ^\ve|\mZ)|  & \leq &
                                \iint\Psi(\mathbf{x}-\mathbf{y})\,\rho_0^\ve(\mathbf{x})\,\left|\rho_0(\mathbf{y})
                                -
                                \rho_0^\ve(\mathbf{y})\right|\,\left|V(t,\mathbf{y})\,-\,V(t,\mathbf{x})\right|\,
                                \left| V(t,\mathbf{x})-V^\ve(t,\mathbf{x})\right| \,\md
                                \mathbf{x}\, \md \mathbf{y}    
\\
                         & \leq&  \|V\|_{L^\infty}
                                 \,\|\Psi\|_{L^1}\,
                                 \left( \|\rho_0^\ve\|_{L^\infty} \|\rho_0 -
                                 \rho_0^\ve\|_{L^2}^2\,+\int\eta
                                 (\mZ^\ve|\mZ)\md \mathbf{x} \right).
\end{eqnarray*}
Gathering these last inequalities, using the uniform control of $\|\rho_0^\ve\|_{L^\infty}$ given by hypothesis \eqref{H2-f0:L1}, we have shown that there exists a
constant $C_T>0$, which does not depend on $\ve>0$ such that for all
$t\in [0,T]$,
\begin{equation}
\label{step:1}
\int_0^t |\mR(\mZ^\ve(s)|\mZ(s))| \md s \,\leq\, C_T \left[ \|\rho_0 -
                                 \rho_0^\ve\|_{L^2}^2 \,+\,
\int_0^t\int\eta(\mZ^\ve(s)|\mZ(s))\,\md \mathbf{x}\,\md s\right].
\end{equation}
It remains to estimate the error term 
\begin{eqnarray*}
\left|\int \left(V^\ve(t)-V(t)\right)\mathcal{E}(f^\ve)(t)\md
  \mathbf{x}\right|  
&\leq&  \frac 3 2\int
\left|V^\ve(t)-V(t)\right|
\left|V^\varepsilon(t) -v \right|\left[ (V^{\varepsilon}(t))^2 +v^2\right]f^{\varepsilon}(t)\md v\md w\md \mathbf{x}
\\
&\leq&\alpha(t) \,\left(\int \left| V^\varepsilon(t)- v\right|^2\,
       f^{\varepsilon}(t)\, \md v\,\md w\,\md \mathbf{x}\right)^{1/2},
\end{eqnarray*}
where $\alpha(t)$ is given by 
$$
\alpha(t)\,:=\, \frac{3}{2}\left(\int  \left[\left(V^\varepsilon(t)\right)^2
  \,+\, v^2\right]^2 \, \left[V^\varepsilon(t)- V(t)\right]^2   \,f^{\varepsilon}(t)\, \md v\,\md w\,\md \mathbf{x}\right)^{1/2}.
$$
Using that $V$ is uniformly bounded in $L^\infty$ according to
Proposition \ref{prop:exhydro} and since 
$$
\rho^\ve_0 \, |V^\ve(t,\mathbf{x})|^6 \leq \int |v|^6 f^\ve(t,\mathbf{x},v,w) \md v \md w
$$
we deduce from Corollary \ref{cor:M2} that there exists a constant $C_T>0$, which does not depend
on $\ve$, such that
$$
\int_0^T\alpha^2(s)\,\md s \leq C_T.
$$
It yields from Lemma \ref{estimate1sansrho} and the Cauchy-Schwarz
inequality that 
\begin{equation}
\label{step:2}
\int_0^T \left|\int \left(V^\ve(t)-V(t)\right)\mathcal{E}(f^\ve)(t)\md
  \mathbf{x}\right|  \md t \leq C_T \,\varepsilon^{{1}/{(d+6)}}.
\end{equation}
Finally integrating (\ref{dd}) on the time interval $[0,t]$ we 
get from the previous estimates  (\ref{step:1}) and (\ref{step:2}) and using the
Gr\"onwall's lemma that
\begin{eqnarray*}
\int \eta({\mathcal{Z}}^\ve(t)|\mathcal{Z}(t))\,\md \mathbf{x} &\leq&  C_T\left[\int
\eta({\mathcal{Z}}^\ve(0)|\mathcal{Z}(0))\,\md \mathbf{x} + \|\rho_0 -
                                 \rho_0^\ve\|_{L^2}^2 +
                                                             \varepsilon^{{1}/{(d+6)}}
                                                             \right]
\end{eqnarray*}
From the assumption \eqref{H5}, we get that for all $t\in [0,T]$,
$$
\int \eta({\mathcal{Z}}^\ve(t)|\mathcal{Z}(t))\,\md \mathbf{x}  \,\leq\, C_T \, \ve^{1/(d+6)},
$$
which concludes the proof.
\end{Proof}

\subsection{Conclusion -- Proof of Theorem \ref{thm:2}}

In this section, we complete the proof of Theorem \ref{thm:2} using
the entropy estimates previously established to show the convergence
of $f^{\varepsilon}$ in the limit $\varepsilon\rightarrow 0$.

First, we set 
$$
F^\ve(t,\mathbf{x},w) \,:=\, \int f^\ve(t,\mathbf{x},v,w) \,\md v,
$$
with an initial datum $F^\ve_0$  given by
$$
F^\ve_0 = \int_{\mathbb{R}} f^\ve _0\,\md v.
$$
Noticing that since  $f^\ve$ is compactly supported in $v$ for any
$\ve >0$, we can choose a test function  in \eqref{weak:sol} independent of
$v\in\mathbb{R}$, hence  the distribution $F^\ve$ satisfies the following equation,
\eqref{eq:feps}
\begin{eqnarray*}
&&\int_0^T\int_{\mathbb{R}^{d+1}} F^\ve \,\partial_t \varphi \,+\, \tau\, \left[\int_{\mathbb{R}} v
  f^\ve \md v \,+\, (a -b\,w) \, F^\ve \right]\,\partial_w\varphi \,\md \mathbf{x}\md w\md t \\
&&+ \int_{\mathbb{R}^{d+1}} F^\ve_0 \,\varphi(0) \,\md \mathbf{x}\md w \,=\, 0, \quad \forall \varphi\in \mathscr{C}_c^\infty([0,T)\times\mathbb{R}^{d+1}),
\end{eqnarray*}
or after reordering
\begin{eqnarray*}
&&\int_0^T\int_{\mathbb{R}^{d+1}} F^\ve \,\left[\partial_t \varphi \,+\, A(V(t,\mathbf{x}),w)\,\partial_w\varphi\right] \,\md \mathbf{x}\md w\md t \,+\, \int_{\mathbb{R}^{d+1}} F^\ve_0 \,\varphi(0) \,\md \mathbf{x}\md w
\\
&& = \tau\int_0^T\int_{\mathbb{R}^{d+2}} \left( V(t,\mathbf{x}) - v\right) f^\ve \,\partial_w\varphi \, \md v \md \mathbf{x}\md w\md t, \quad \forall \varphi\in \mathscr{C}_c^1([0,T)\times\mathbb{R}^{d+1}),
\end{eqnarray*}
where $V$ is solution to \eqref{eq:lim}.

On the one hand, using that up to a subsequence $F^\ve$ converges
weakly-$\star$ in $\mathcal{M}((0,T)\times\mathbb{R}^{d+1})$ to a
limit $F \in \mathcal{M}((0,T)\times\mathbb{R}^{d+1})$, we can pass to the limit on the left hand side by
linearity. On the other hand,  from Lemma 
\ref{corD1sansrho} and Proposition \ref{prop:estimentrop},  we get when
$\ve\rightarrow 0$,
\begin{eqnarray*}
\int_0^T \int f^{\varepsilon} |v-V(t,\mathbf{x})|^2 \md \mathbf{x} \md v \md w \md s &
                                                                   \leq&
                                                                         \int_0^T \int f^{\varepsilon} \left(|v-V^{\varepsilon}(t,\mathbf{x})|^2 + |V^{\varepsilon}(t,\mathbf{x})-V(t,\mathbf{x})|^2 \right) \md \mathbf{z} \md s 
\\
& \leq& C_T \,\ve^{1/(d+6)},
\end{eqnarray*}
hence it yields that since $\rho^\ve$ does not depend on time,
$$
\left|\int \left (V(t,\mathbf{x}) - v\right) f^\ve \,\partial_w\varphi \, \md v \md
\mathbf{x}\md w\md t\right| \,\leq C_T\,\|\partial_w\varphi\|_{L^\infty}\,\|\rho_0^\ve\|_{L^1}^{1/2} \, \ve^{1/(2d+12)}.
$$
Thus, passing to the limit $\ve\rightarrow 0$, it proves that $F$ is a measure solution
\eqref{eq:F}. Furthermore, by uniqueness of the solution to
\eqref{eq:F}, we get the convergence for the sequence
$(F^\ve)_{\ve>0}$.

Finally let us show that for any $\varphi\in \mathscr{C}^0_b(\R^{d+2})$, 
$$ 
\ds \int \varphi(\mathbf{x},v,w)\, f^{\varepsilon}(t,\mathbf{x},v,w)\,\md v\, \md w\, \md \mathbf{x} \rightarrow \int \varphi(\mathbf{x},V(t,\mathbf{x}),w)\, F(t,\md \mathbf{x},\md w)\, ,  
$$
strongly in
$L^1_{\text{loc}}(0,T)$ as $\varepsilon \rightarrow 0$. Consider $0<t<t'\leq T$. We start with showing the convergence for any $\varphi \in \mathscr{C}^1_c(\R^{d+2})$, and then we will conclude using a density argument. Consider $\varphi \in \mathscr{C}^1_c(\R^{d+2})$, we have:
\begin{align*}
\mathcal{I} \, &:= \, \ds \int_t^{t'}\left|\int \varphi(\mathbf{x},v,w)\, f^{\varepsilon}(s,\mathbf{x},v,w)\,\md v\, \md w\, \md \mathbf{x}  \,-\, \int \varphi(\mathbf{x},V(s,\mathbf{x}),w)\, F(s,\md \mathbf{x},\md w)  \right|\,\md s   \\
& \leq \mathcal{I}_1 \,+\, \mathcal{I}_2  , \\
\end{align*}
where 
$$\left\{
\begin{array}{l l l}
\mathcal{I}_1 &:=& \ds \int_t^{t'}\left|\int \left(\varphi(\mathbf{x},v,w) \,-\,\varphi(\mathbf{x},V(s,\mathbf{x}),w) \right)\, f^{\varepsilon}(s,\mathbf{x},v,w)\,\md v\, \md w\, \md \mathbf{x}   \right|\,\md s  , \\
\mathcal{I}_2 &:=& \ds \int_t^{t'}\left|\int \varphi(\mathbf{x},V(s,\mathbf{x}),w) \, \left(F^{\varepsilon}(s,\md \mathbf{x},\md w)  \,-\, F(s,\md \mathbf{x},\md w) \right)   \right|\,\md s  . \\
\end{array}
\right.$$
We estimate the first term $\mathcal{I}_1$ using the regularity of $\varphi$ and the Cauchy-Schwarz inequality:
\begin{align*}
\mathcal{I}_1 \,&\leq \, \ds \|\partial_v \varphi\|_{L^\infty} \int_0^{T} \int \left| v\,-\,V(s,\mathbf{x}) \right|\,f^\ve(s,\mathbf{x},v,w)\, \md v\, \md w\, \md \mathbf{x} \, \md s   \\
& \leq \,\ds \|\partial_v \varphi\|_{L^\infty} \,\left( \int_0^T\int f^\ve(s,\mathbf{x},v,w)\, \md v\, \md w\, \md \mathbf{x} \, \md s \right)^{1/2}   \\
& ~~~\,\times\,\left( \int_0^T\int f^\ve(s,\mathbf{x},v,w)\,\left| v\,-\,V(s,\mathbf{x}) \right|^2\, \md v\, \md w\, \md \mathbf{x} \, \md s \right)^{1/2}   \\
& \leq \|\partial_v \varphi\|_{L^\infty} \, \left(T \, \|f_0^\ve\|_{L^1} \,C_T\right)^{1/2}\,\ve^{1/(2d+12)},
\end{align*}
whereas the second term $\mathcal{I}_2$ also converges to zero when $\ve$ goes to zero since $F^\ve$ converges weakly-$\star$ in $\mathcal{M}((0,T)\times\mathbb{R}^{d+1})$ to $F$. Using a density argument, this shows the convergence of $f^{\varepsilon}$ in $L^1_{\text{loc}}((0,T),\mathcal{M}(\mathbb{R}^{d+2}))$, so this concludes the proof of Theorem \ref{thm:2}.


%% file: AppThm1.tex
This appendix is devoted to the proof of the existence and uniqueness of a solution $f^\ve$ to \eqref{eq:feps}. Let $T>0$ and $\ve>0$ be fixed. The main difficulty is that we cannot use a compactness argument based on an average lemma as in \cite{KARPER1} for example, since there is no transport term of the form $v\,\nabla_x f^\ve$ in \eqref{eq:feps}. Thus, our strategy is to linearize the equation \eqref{eq:feps} in order to construct a Cauchy sequence which converges towards a solution to \eqref{eq:feps}. First of all, we need to prove some extra a priori estimates on $f^\ve$. In the rest of this section, since $\ve$ is fixed, for the sake of clarity, we will note $f$ instead of $f^\ve$, but all the following estimates are not uniform in $\ve$.

\subsection{A priori estimates}

We recall that for any solution $f$ of \eqref{eq:feps} and for all $t\in[0,T]$, $\|f(t)\|_{L^1} = \|f(0)\|_{L^1}$. This subsection is devoted to the proof of some a priori estimates. 

First of all, we set the system of characteristic equations associated to \eqref{eq:feps} for all $(\mathbf{x},v,w)\in\R^{d+2}$ and all $t\in[0,T]$:
\begin{equation}\label{char}
\left\{\begin{array}{l}
\dfrac{\md \Sigma^v(s)}{\md s}  \,=\, N(\Sigma^v(s)) \,-\, \Sigma^w(s) \,-\,[\Phi_\ve\ast\rho_0](\mathbf{x}) \Sigma^v(s) \,+\, [\Phi_\ve\ast j_f(s)](\mathbf{x}),   \\ \, \\
\dfrac{\md \Sigma^w(s)}{\md s}  \,=\, A(\Sigma^v(s),\Sigma^w(s)),   \\ \, \\
\Sigma^v(t) \,=\, v , \quad \Sigma^w(t) \,=\, w,
\end{array}  
\right.
\end{equation}
where $f$ is a solution of \eqref{eq:feps} and $j_f$ is defined with:
$$j_f:=\ds\int f\,v\,\md v\,\md w.$$
We define the flow of \eqref{char} for all $\mathbf{z}=(\mathbf{x},v,w)\in\R^{d+2}$ for all $s,t$ in $[0,T]$:
$$\Sigma(s,t,\mathbf{z})=\left( \Sigma^v(s,t,\mathbf{z}) \, , \, \Sigma^w(s,t,\mathbf{z}) \right) \quad , \quad \Sigma(t,t,\mathbf{z})=\left( v \, , \, w \right).$$
We start with proving the well-posedness of the flow of the characteristic equation \eqref{char} and an estimate of the support of a solution to \eqref{eq:feps}.
\begin{lem}[Well-posedness of the characteristic system and estimate of the support]\label{estimsupp}
Consider an initial data $f_0$ satisfying \eqref{H3-1} and \eqref{H3-2}, and suppose that there exists $f$ a smooth non-negative solution to \eqref{eq:feps} such that for all $t\in[0;T]$, $\|f(t)\|_{L^1} = \|f(0)\|_{L^1}$. Then, the characteristic system \eqref{char} is well-posed. Furthermore, there exists a positive constant $R_{T,\ve}$ such that:
\begin{equation}\label{supp}
\underset{\mathbf{x}\in\R^d}{\underset{t\in[0,T]}{\sup}} \Supp(f(t,\mathbf{x},\cdot,\cdot))  \subset B(0,R_{T,\ve}),
\end{equation}
and there exist two positive constants $C_1$ and $C_2$ such that:
\begin{equation}\label{suppunif}
R_{T,\ve}\leq C_2 e^{C_1(1+\frac{1}{\ve})T} .
\end{equation}
\end{lem}

\begin{Proof}
The Cauchy-Lipschitz Theorem yields the local existence and uniqueness of the flow of the characteristic equation \eqref{char}. Define for all $s\in[0,T]$:
$$R_\ve(s) := \sup\left\{\|\Sigma(s',0,\mathbf{z})\| \quad|\quad  s'\in[0;s],\, \mathbf{z}\in\R^d\times B(0,R_0^\ve) \right\}.$$
Our purpose is to estimate $R_\ve$ using the following energy estimate. Let $s\in[0;T]$, $\mathbf{z}=(\mathbf{x},v,w)\in\R^d\times B(0,R_0^\ve)$. We have:

\begin{align*}
\|\Sigma(s,0,\mathbf{z})\|^2 \, & =\, \|(v,w)\|^2 \,+\, 2\displaystyle\int_0^s \Sigma(s',0,\mathbf{z}) \, \partial_s \Sigma(s',0,\mathbf{z}) \, \md s' \\
& =\, \|(v,w)\|^2 \,+\, 2\displaystyle\int_0^s \left[\Sigma^v \, N(\Sigma^v) - \Sigma^v \, \Sigma^w  \right. \\
& ~~ \left. -\, [\Phi_\ve\ast\rho_0](\mathbf{x})\,|\Sigma^v|^2 \,+\, [\Phi_\ve\ast j_f(s')](\mathbf{x})\,\Sigma^v     \,+\, \Sigma^w\, A(\Sigma^v,\Sigma^w) \right] \md s'   \\
& \leq\, |R_0^\ve|^2 \,+\, 2\displaystyle\int_0^s \left[ |\Sigma^v|^2 + \dfrac{1}{2}\, (|\Sigma^v|^2 + |\Sigma^w|^2)  \right. \\
& ~~ \left. +\, R_\ve(s)\,[\phi_\ve\ast\rho_0](\mathbf{x})\,\Sigma^v  \,+\, \dfrac{\tau}{2}\,(|\Sigma^v|^2 + |\Sigma^w|^2) \,+\, \dfrac{\tau}{2}\,(a^2 + |\Sigma^w|^2) \right] \md s'  \\
& \leq\, |R_0^\ve|^2 \,+\, T\,\tau\, a^2 \,+\, \displaystyle\int_0^s\left( 3 + 2\,\tau + 2\,\|\rho_0\|_{L^\infty}\,\|\Psi\|_{L^{1}} + \frac{2}{\ve} \, \|\rho_0\|_{L^\infty} \right)\,|R_\ve(s')|^2 \, \md s'.
\end{align*}
Thus, passing to the supremum in $\mathbf{x}\in\R^d$ and in $(v,w)\in B(0,R_0^\ve)$, we get that for all $s\in[0,T]$:
$$|R_\ve(s)|^2 \,\leq\, |R_0^\ve|^2 \,+\, T\,\tau\, a^2 \,+\, \displaystyle\int_0^s \left( 3 + 2\, \tau + 2\, \|\rho_0\|_{L^\infty}\, \|\Psi\|_{L^{1}} + \frac{2}{\ve} \, \|\rho_0\|_{L^\infty} \right)R_\ve(s')\, \md s'.$$
Using the Gr\"owall's inequality, we get that there exist two positive constants $C_1$ and $C_2$ which depend only on $T$, $R_0^\ve$, $\|\Psi\|_{L^1}$, $\rho_0$, $\tau$ and $a$ such that for all $s\in[0,T]$:
$$|R_\ve(s)|^2 \,\leq\, C_2 \, e^{C_1 \,(1+\frac{1}{\ve})\,T} .$$
Therefore, the function $\Sigma(\cdot,0,\cdot)$ is well-posed in $\mathscr{C}^1([0,T]^2\times B(0,R_0^\ve) \,,\, L^\infty(\R^d_\mathbf{x}))$. Similarly, for all $t\in[0,T]$, we show that $\Sigma(\cdot,t,\cdot)$ is well-posed in $\mathscr{C}^1([0,T]^2\times B(0,R_\ve(t)) \,,\, L^\infty(\R^d_\mathbf{x}))$. Furthermore, we can conclude that
$$\underset{\mathbf{x}\in\R^d}{\underset{t\in[0,T]}{\sup}} \Supp(f(t,\mathbf{x},\cdot,\cdot))  \subset B(0,R_\ve(T)).$$
\end{Proof}

Now, using this estimate on the propagation of the support in $\mathbf{u}=(v,w)$ of any solution $f$ to \eqref{eq:feps}, we can estimate $f$ and $\nabla_\mathbf{u} f$ in $L^\infty$.

\begin{lem}[Estimates in $L^\infty$]\label{lem:Linf}
Consider an initial data $f_0$ satisfying \eqref{H3-1} and \eqref{H3-2}, and suppose that there exists $f$ a smooth non-negative solution to \eqref{eq:feps} such that for all $t\in[0,T]$, $\|f(t)\|_{L^1} = \|f_0\|_{L^1}$. Then, there exists a positive constant $C_{T,\ve}$ such that for all $t\in[0,T]$,
$$\|f(t)\|_{L^{\infty}} \leq C_{T,\ve} \quad \text{, and }\quad \|\nabla_\mathbf{u} f(t)\|_{L^{\infty}}\leq C_{T,\ve}.$$ 
\end{lem}

\begin{Proof}
We write \eqref{eq:feps} in a non-conservative form: 
\begin{equation}\label{noncons}
\partial_t f + \mathbf{A}\, \cdot\, \nabla_\mathbf{u} f = -\displaystyle\dv\ds_\mathbf{u} \left( \mathbf{A} \right)f,
\end{equation}
where $\mathbf{A}$ is the advection field of \eqref{eq:feps} given for all $t\in[0,T]$ and all $\mathbf{z}=(\mathbf{x},v,w)\in\R^{d+2}$ by:
$$\mathbf{A}(t,\mathbf{z}):=\begin{pmatrix} N(v) \,-\, w \,-\,[\Phi_\ve\ast\rho_0](\mathbf{x})\, v \,+\, [\Phi_\ve\ast j_f(t)](\mathbf{x}) \\ \, \\ A(v,w) \end{pmatrix}.$$
Thus, for all $t\in[0,T]$ and all $\mathbf{z}=(\mathbf{x},v,w)\in\R^{d+2}$, we note $\mathbf{u}=(v,w)$ and we have:
\begin{align*}
 -\dv\ds_\mathbf{u} \left( \mathbf{A}(t,\mathbf{z}) \right)\,& = \, -N'(v) \,+\, [\Phi_\ve\ast\rho_0](\mathbf{x}) \,+\, \tau\,b   \\
 & = \, 3\,v^2 \,-\,1 \,+\, [\Phi_\ve\ast\rho_0](\mathbf{x}) \,+\, \tau\,b .
\end{align*}
Then, we get that for all $t\in[0,T]$ and $\mathbf{z}=(\mathbf{x},v,w)\in\R^{d+2}$:
$$f(t,\mathbf{z}) = f_0(\mathbf{x},\Sigma(0,t,\mathbf{z})) \,-\, \displaystyle\int_0^t \left( \dv\ds_\mathbf{u} \left( \mathbf{A} \right) \, f \right)(t,\mathbf{x},\Sigma(s,t,\mathbf{z})) \, \md s .$$
Consequently, we have for all $t\in[0,T]$: 
\begin{align*}
\|f(t)\|_{L^\infty}\, &\leq\, \|f_0\|_{L^\infty} + \displaystyle\int_0^t \|\dv\ds_\mathbf{u}  \left( \mathbf{A}(s,\cdot) \right) f(s)\|_{L^\infty}\, \md s      ,
\end{align*}
Moreover, according to Lemma \ref{estimsupp}, there exists a positive constant $R_{T,\ve}$ satisfying \eqref{supp}, and therefore, for all $s$, $t\in[0,T]$ and $\mathbf{z}=(\mathbf{x},\mathbf{u})\in\R^d\times B(0,R_{T,\ve})$,
$$|\dv\ds_\mathbf{u}  \left( \mathbf{A}(s,\mathbf{z}) \right)|  \,\leq \,3\,|R_{T,\ve}|^2 \,+\,1 \,+\, \|\Phi_\ve\|_{L^1}\,\|\rho_0\|_{L^\infty} \,+\,\tau\,b  .$$
Therefore, the Gr\"onwall's inequality gives us that there exists a positive constant $C_1$ such that for all $t\in[0,T]$:
$$\|f(t)\|_{L^\infty}\,\leq\, \|f_0\|_{L^\infty} e^{ C_1 \, t \, } .$$
Then, by differentiating \eqref{noncons} with respect to $v$ and $w$, we have:
$$
\left\{
\begin{array}{l l l}
\partial_t\left(\partial_v f\right)  + \mathbf{A}\, \cdot\, \nabla_\mathbf{u} \left(\partial_v f\right) &=& S^v(t,\mathbf{z})\,\partial_v f - \tau \,\partial_w f - N''(v)\,f ,   \\ \, \\
\partial_t\left(\partial_w f\right)  + \mathbf{A}\, \cdot\, \nabla_\mathbf{u} \left(\partial_w f\right) &=&  S^w(t,\mathbf{z})\,\partial_w f + \partial_v f ,
\end{array}
\right.$$
where $S^v$ and $S^w$ are given for all $t\in[0,T]$ and all $\mathbf{z}=(\mathbf{x},v,w)\in\R^{d+2}$ by:
$$\left\{
\begin{array}{l}
S^v(t,\mathbf{z})\,:=\, -2\, N'(v)  + 2\,\left[\Phi_\ve\ast\rho_0\right](\mathbf{x}) + \tau\,b ,\\ \, \\
S^w(t,\mathbf{z})\,:=\, -N'(v) + \left[\Phi_\ve\ast\rho_0\right](\mathbf{x}) + 2\,\tau\,b.
\end{array}
\right.$$
Therefore, using the estimate of the support of $f$ \eqref{supp}, the Gr\"onwall's inequality gives us that there exist two positive constants $C_2$ and $C_3$ such that for all $t\in[0,T]$:
$$
\|\nabla_\mathbf{u} f(t)\|_{L^\infty} \,\leq\, C_3\, e^{C_2\, t}   .
$$

\end{Proof}

\subsection{Proof of existence and uniqueness}

We proceed with a linearization of the equation \eqref{eq:feps} in order to construct a Cauchy sequence $\{f^\ve_n\}_{n\in\mathbb{N}}$ of non-negative functions in the Banach space $\mathscr{C}^0([0,T]\times \R^2_\mathbf{u}\,,\,L^\infty(\R^d_\mathbf{x}))$, such that for all $n\in\mathbb{N}$, the Lemmas \ref{estimsupp} and \ref{lem:Linf} give that there exist two positive constants $R_{T,\ve}$ and $C_{T,\ve}$ independent on $n$ such that for all $n\in\mathbb{N}$:
$$\underset{\mathbf{x}\in\R^d}{\underset{t\in[0,T]}{\sup}} \Supp(f^\ve_n(t,\mathbf{x},\cdot,\cdot))  \subset B(0,R_{T,\ve}),$$
$$\|f^\ve_n(t)\|_{L^{\infty}} \leq C_{T,\ve} \quad \text{ and }\quad \|\nabla_\mathbf{u} f^\ve_n(t)\|_{L^{\infty}}\leq C_{T,\ve}.$$ 
Thus by a classical fixed point argument in a Banach space, the sequence $\{f^\ve_n\}_{n\in\mathbb{N}}$ strongly converges in $\mathscr{C}^0([0,T]\times \R^2_\mathbf{u}\,,\,L^\infty(\R^d_\mathbf{x}))$ towards a function $f^\ve$, which is non-negative and
$$ f^\ve\in\mathscr{C}^0([0,T]\times \R^2_\mathbf{u}\,,\,L^\infty(\R^d_\mathbf{x})).$$
Hence, passing to the limit $n\rightarrow+\infty$ in the linearized equation, $f^\ve$ is a weak solution to \eqref{eq:feps} in the sense of (\ref{weak:sol}) which satisfies the support estimate \eqref{supp}. Consequently, we deduce
$$f^\ve\in \mathscr{C}^0\left( [0,T] \,,\, L^1(\R^{d+2}) \right),$$
hence we can apply Lemma \ref{estimsupp} with $f^\ve$, which yields that the characteristics $\Sigma \in \mathscr{C}^1([0,T]^2\times \R^2_\mathbf{u} \,,\, L^\infty(\R^d_\mathbf{x}))$ are well-defined. As previously, we note $\mathbf{A}$ the advection field of \eqref{eq:feps}. Therefore, for all $t\in[0,T]$, $\mathbf{u}=(v,w)\in \R^2$ and $\mathbf{x}\in\R^d$,
$$f^\ve(t,\mathbf{x},\mathbf{u}) \,=\, f_0^\ve(\mathbf{x},\Sigma(0,t,\mathbf{x},\mathbf{u})) \, \exp\left(-\ds\int_0^t \dv\ds_\mathbf{u} \mathbf{A}(s,\mathbf{x},\Sigma(s,t,\mathbf{x},\mathbf{u})) \md s  \right) .$$
Finally, using the regularity of $\dv\ds_\mathbf{u} \mathbf{A}$, $\Sigma$ and $f_0^\ve$ which satisfies \eqref{H3-1}, we get that
$$f^\ve, \, \nabla_{\mathbf{u}} f^\ve \in L^\infty((0,T)\times\mathbb{R}^{d+2}).$$

%% file: AppThm2.tex
This appendix is devoted to the proof of Proposition
\ref{prop:exhydro}. We apply a fixed point argument to get the
existence and uniqueness of a classical solution $(V,W)$ of
\eqref{eq:nlFHN}. We first  set 
$$
\mathscr{E}\,:=\,\mathscr{C}^0([0,T],L^{\infty}(\R^d))
$$ 
and for a fixed $M>0$ we define
$$
\mathscr{K}_{T}\,:=\,\left \{ (V,W)\in \mathscr{E}^2 \text{ ; } \|V(t) -
  V_0\|_{L^{\infty}(\mathbb{R}^d)}+\|W(t) -
  W_0\|_{L^{\infty}(\mathbb{R}^d)} \leq M, \, \forall
  t\in[0,T] \right \},
$$
equipped with the norm for $U=(V,W)\in \mathscr{K}_{T}$,
\begin{equation*}
\vertiii{ U} \, :=\,\underset{t \in [0,T]}{\sup} \left( \|V(t)\|_{L^{\infty}(\mathbb{R}^d)}+\|W(t) \|_{L^{\infty}(\mathbb{R}^d)}\right).
\end{equation*}
Consider the application $\Gamma$ such that for all $U=(V,W)\in \mathscr{K}_T$,
for all $t\in[0,T]$ and almost every $\mx\in\mathbb{R}^d$,
$$
\Gamma[U](t,\mx)\,:=\,U_0(\mx) \,+\, \displaystyle\int_0^t 
\left[
\begin{array}{l}
\ds\mathcal{L}_{\rho_0}(V)(s,\mx)\md s + N(V(s,\mx)) - W(s,\mx) 
\\ \, \\
A(V,W)(s,\mx)
\end{array}
\right] \md s.$$
On the one hand, since $\Psi\in L^1(\mathbb{R}^d)$ and $N \in
\mathscr{C}^1(\mathbb{R})$,  we   prove that $\Gamma$ is well-defined
$$
\vertiii{\Gamma[U] - U_0} \, \leq \, C\, T 
$$
for some constant $C>0$ independent of time. As a consequence, we can
choose $T$ small enough so that $\Gamma[U] \in \mathscr{K}_{T}$.  On the other
hand, $\Gamma$ is contractive : for $U_1$  $U_2\in \mathscr{K}_{T}$ and $(t,\mx)$ be
in $[0,T]\times\mathbb{R}^d$, we have
\begin{equation*}
\left|\Gamma[U_1] \,-\, \Gamma[U_2]\right|(t,\mx)   \leq C \,t\, \|U_1-U_2\|,
\end{equation*}
where $C>0$ does not depend of time $t>0$, hence for $T$ small enough,
$\Gamma$ is a contraction from $\mathscr{K}_T$ to $\mathscr{K}_T$, which is a complete
for the distance associated to the norm $\vertiii{.}$ previously
defined. Therefore, using the Banach fixed point theorem, there exists
a unique $U\in \mathscr{K}_T$ such that $\Gamma[U]=U$, that is,  $U=(V,W)$ is
solution to (\ref{eq:nlFHN}) for the initial condition
$U_0=(V_0,W_0)$. Moreover since $V$ and $W\in
\mathscr{C}^0([0,T],L^{\infty}(\mathbb{R}^d))$ and using  the
Duhamel's  formula, it yields that 
$$
V\, , \,W\in \mathscr{C}^1([0,T], L^{\infty}(\mathbb{R}^d)).
$$ 
Finally, an energy estimate gives that
\begin{eqnarray*}
\dfrac{1}{2}\,\frac{\md }{\md t}\left( |V(t,\mx)|^2+  |W(t,\mx)|^2\right)
  &\leq&  \left(\dfrac{3+\tau}{2} +
         \|\Psi\|_{L^1}\,\|\rho_0\|_{L^\infty}
         \right)\|V(t)\|_{L^{\infty}(\mathbb{R}^d)}^2 
\\ &+&
         \frac{1+2\,\tau}{2}\, \|W(t)\|_{L^{\infty}(\mathbb{R}^d)}^2  +
         \frac{\tau \,a^2}{2}.
\end{eqnarray*}
Using Gr\"onwall's lemma, we can conclude that there exists a constant
$C>0$, only depending on $a$, $\tau$ and $\Psi$, such that
\begin{equation*}
\vertiii{(V,W)} \,\leq\,   \left( \vertiii{(V_0,W_0)} + C\right)\, e^{C\,T}.
\end{equation*}
It  allows to prove that this unique solution
is global in time.